\documentclass[12pt]{article}

\usepackage{amssymb}
\usepackage{luis}

\def\stackunder#1#2{\mathrel{\mathop{#2}\limits_{#1}}}
\def\QTR#1#2{{\csname#1\endcsname #2}}

\def\dfrac#1#2{{\displaystyle {#1 \over #2}}}
\def\QATOP#1#2{{#1 \atop #2}}
\def\tfrac#1#2{{\textstyle {#1 \over #2}}}

\newtheorem{theorem}{Theorem}[section]

\newtheorem{lemma}[theorem]{Lemma}
\newtheorem{proposition}[theorem]{Proposition}
\newtheorem{definition}[theorem]{Definition}
\newtheorem{example}[theorem]{Example}

\newtheorem{remark}[theorem]{Remark}

\begin{document}

\vspace*{3cm}
\begin{center}
{\LARGE Analysis on Poisson and Gamma spaces\bigskip }

{\sc Yuri G.~Kondratiev}$^{1,2,3}$

{\sc Jos\'{e} L.~Silva}$^{2,4}$

{\sc Ludwig Streit}$^{2,4}$

{\sc Georgi F.~Us}$^{5}$
\bigskip 

\begin{tabular}{l}
$^{1}$Inst.~Angewandte Math., Bonn Univ., D-53115 Bonn, Germany\\
$^{2}$BiBoS, Bielefeld Univ., D-33615 Bielefeld, Germany \\ 
$^{3}$Inst.~Math., 252601 Kiev, Ukraine\\ 
$^{4}$CCM, Univ.~Madeira, P-9000 Funchal, Portugal \\
$^{5}$Dep.~ Mech.~ and Math.~Univ.~Kiev, Kiev, 252033 Ukraine \\
\multicolumn{1}{r}{CCM$^{@}$-UMa 24/97}
\end{tabular}
\end{center}

\renewcommand{\thefootnote}{}
\footnotetext{
$^{@}$http:/www.uma.pt/ccm/ccm/ccm.html}

\thispagestyle{empty}

\newpage
\renewcommand{\thefootnote}{\arabic{footnote}}
\setcounter{page}{1}

\title{Analysis on Poisson and Gamma spaces}
\author{\textbf{Yuri G. Kondratiev} \\
BiBoS, Universit\"at Bielefeld, D 33615 Bielefeld, Germany \\
Institute of Mathematics, Kiev, Ukraine 
\and 
\textbf{Jos\'{e} L. da Silva} \\
CCM, Universidade da Madeira, P 9000 Funchal, Portugal\\
(luis@uma.pt)
\and 
\textbf{Ludwig Streit} \\
BiBoS, Universit\"at Bielefeld, D 33615 Bielefeld, Germany \\
CCM, Universidade da Madeira, P 9000 Funchal, Portugal 
\and 
\textbf{Georgi F. Us} \\
Department of Mechanics and Mathematics, University Kiev, \\
Kiev, 252033 Ukraine}
\date{}
\maketitle

\begin{abstract}
We study the spaces of Poisson, compound Poisson and Gamma noises as special
cases of a general approach to non-Gaussian white noise calculus, see \cite
{KSS96}. We use a known unitary isomorphism between Poisson and compound
Poisson spaces in order to transport analytic structures from Poisson space
to compound Poisson space. Finally we study a Fock type structure of chaos
decomposition on Gamma space.
\end{abstract}

\renewcommand{\thefootnote}{}
\footnotetext{
Published in \textit{Infinite Dimensional Analysis, Quantum Probabilities
and Related Topics}. Vol.~{\bf 1},N.1, pp.~91-117, 1998.}

\newpage
\renewcommand{\thefootnote}{\arabic{footnote}}

\tableofcontents

\section{Introduction}

The present paper elaborate the $L^{2}$ structure of compound Poisson
spaces; we note that for all of this compound Poisson processes the results
of \cite{KSS96} immediately produce Gel'fand triples of test and generalized
functions as well their characterizations and calculus.

The Analysis on pure Poisson spaces was developed in \nocite{ABT90} \cite
{CP90}, \cite{IK88}, \cite{NV95}, \cite{P95}, \cite{U95} and many others
from different points of view. In \cite{KSS96} we have developed methods for
non Gaussian analysis based on generalized Appell systems. In the case of
Poisson space, this coincide with the system of generalized Charlier
polynomials, however the desirable extensions to compound Poisson and for
example Gamma processes are trivial.

Let us describe this construction more precisely. We recall that the Poisson
measure $\pi _{\sigma }$ (with intensity measure $\sigma $ which is a
non-atomic Radon measure on $\QTR{mathbb}{R}^{d}$) is defined by its Laplace
transform as 
\[
l_{\pi _{\sigma }}(\varphi )=\int_{\mathcal{D}^{\prime }}\exp \left(
\left\langle \gamma ,\varphi \right\rangle \right) d\pi _{\sigma }(\gamma
)=\exp \left( \int_{\QTR{mathbb}{R}^{d}}(e^{\varphi (x)}-1)d\sigma
(x)\right) ,\,\varphi \in \mathcal{D}, 
\]
where $\mathcal{D}^{\prime }$ is the dual of $\mathcal{D}:=\mathcal{D}(%
\QTR{mathbb}{R}^{d})=C_{0}^{\infty }(\QTR{mathbb}{R}^{d})$ ($C^{\infty }$%
-functions on $\QTR{mathbb}{R}^{d} $ with compact support). An additional
analysis shows that the support of the measure $\pi _{\sigma }$ consists of
generalized functions of the form $\sum_{x\in \gamma }\varepsilon _{x},$ $%
\gamma \in \Gamma _{\QTR{mathbb}{R}^{d}},$ where $\varepsilon _{x}$ is the
Dirac measure in $x$ and $\Gamma _{\QTR{mathbb}{R}^{d}}$ is the
configuration space over $\QTR{mathbb}{R}^{d},$ i.e., 
\[
\Gamma _{\QTR{mathbb}{R}^{d}}:=\{\gamma \subset \QTR{mathbb}{R}%
^{d}\,|\,\left| \gamma \cap K\right| <\infty \;\mathrm{for\;any\;compact\;}%
K\subset \QTR{mathbb}{R}^{d}\}. 
\]
The configuration space $\Gamma _{\QTR{mathbb}{R}^{d}}$ can be endowed with
its natural Borel $\sigma $-algebra $\mathcal{B}(\Gamma _{\QTR{mathbb}{R}%
^{d}})$ and $\pi _{\sigma }$ can be considered as a measure on $\Gamma _{%
\QTR{mathbb}{R}^{d}}.$

Let us choose a transformation $\alpha $ on $\mathcal{D}$ given by 
\[
\alpha (\varphi )(x)=\log (1+\varphi (x)),\;-1<\varphi \in \mathcal{D}%
,\,x\in \QTR{mathbb}{R}^{d}.
\]
Then the normalized exponential or Poisson exponential 
\[
e_{\pi _{\sigma }}^{\sigma }(\varphi ;\gamma )=\exp \left( \left\langle
\gamma ,\log (1+\varphi )\right\rangle -\left\langle \varphi \right\rangle
_{\sigma }\right) ,\;\gamma \in \Gamma _{\QTR{mathbb}{R}^{d}}
\]
is a real holomorphic function of $\varphi $ on a neighborhood of zero $%
\mathcal{U}_{\alpha }$ on $\mathcal{D}.$ Its Taylor decomposition (with
respect to $\varphi $) has the form 
\[
e_{\pi _{\sigma }}^{\sigma }(\varphi ;\gamma )=\sum_{n=0}^{\infty }\frac{1}{%
n!}\left\langle C_{n}^{\sigma }(\gamma ),\varphi ^{\otimes n}\right\rangle
,\;\varphi \in \mathcal{U}_{\alpha }^{\prime }\subset \mathcal{U}_{\alpha
},\,\gamma \in \Gamma _{\QTR{mathbb}{R}^{d}},
\]
where $C_{n}^{\sigma }:\Gamma _{\QTR{mathbb}{R}^{d}}\rightarrow \mathcal{D}%
^{\prime \widehat{\otimes }n}.$ It follows from the above equality that for
any $\varphi ^{(n)}\in \mathcal{D}^{\widehat{\otimes }n},$ $n\in 
\QTR{mathbb}{N}_{0}$ the function 
\[
\Gamma _{\QTR{mathbb}{R}^{d}}\ni \gamma \mapsto \left\langle C_{n}^{\sigma
}(\gamma ),\varphi ^{(n)}\right\rangle 
\]
is a polynomial of order $n$ on $\Gamma _{\QTR{mathbb}{R}^{d}}.$ This is
precisely the system of generalized Charlier polynomials for the measure $%
\pi _{\sigma },$ see Subsection \ref{2eq61} for details.

This system can be used for a Fock realization. $L^{2}(\pi _{\sigma })$ has
a Fock realization analogous to Gaussian analysis, i.e, 
\[
L^{2}(\pi _{\sigma })\simeq \bigoplus_{n=0}^{\infty }\mathrm{Exp}%
_{n}L^{2}(\sigma )=\mathrm{Exp}L^{2}(\sigma ), 
\]
where $\mathrm{Exp}_{n}L^{2}(\sigma )$ denotes the n-fold symmetric tensor
product of $L^{2}(\sigma ).$

The ``Poissonian gradient'' $\bigtriangledown ^{\mathtt{P}}$ on functions $%
f:\Gamma _{\QTR{mathbb}{R}^{d}}\rightarrow \QTR{mathbb}{R}$ which has
specific useful properties on Poisson space, is introduced on a specific
space of ``nice'' functions as a difference operator 
\[
(\bigtriangledown ^{\mathtt{P}}f)(\gamma ;x)=f\left( x+\varepsilon
_{x}\right) -f(x),\;\gamma \in \Gamma _{\QTR{mathbb}{R}^{d}},\,x\in 
\QTR{mathbb}{R}^{d}. 
\]
The gradient $\bigtriangledown ^{\mathtt{P}}$ appears from different points
of view in many papers on conventional Poissonian analysis, see e.g. \cite
{IK88}, \cite{NV95}, \cite{KSS96} and references therein. We note also that
the most important feature of the Poissonian gradient is that it produces
(via a corresponding integration by parts formula) the orthogonal system of
generalized Charlier polynomials on $(\Gamma _{\QTR{mathbb}{R}^{d}},\mathcal{%
B}(\Gamma _{\QTR{mathbb}{R}^{d}}),$\linebreak $\pi _{\sigma }),$ see Remark 
\ref{2eq62}. In addition we mention here that as a tangent space to each
point $\gamma \in \Gamma _{\QTR{mathbb}{R}^{d}}$ we choose the same Hilbert
space $L^{2}(\QTR{mathbb}{R}^{d},\sigma ).$

We conclude Section \ref{ref21} with the expressions for the annihilation
and creation operators on Poisson space. In terms of chaos decomposition
Nualart and Vives \cite{NV95} proved the analogous expression of the
creation operator, in this paper we give an independent proof which is based
on the results on absolute continuity of Poisson measures, see e.g. \cite
{Sk57} and \cite{T90}, details can be found in Subsection \ref{2eq63}.

The analysis on compound Poisson space can be done with the help of the
analysis derived from Poisson space described above. That possibility is
based on the existence of an unitary isomorphism between compound Poisson
space and Poisson space which allows us to transport the Fock structure from
the Poisson space to the compound Poisson space. The above isomorphism has
been identified before by K. It\^{o}, \cite{I56} and A. Dermoune, \cite{De90}%
. All this is developed is Subsection \ref{2eq64}.

The images of the annihilation and creation operators under the above
isomorphism on compound Poisson space are worked out in Subsection \ref{2eq65}%
.

The aim of Section \ref{2eq66} is to study in more details the previous
analysis in a particular case of compound Poisson measure, the so called
Gamma noise measure. Its Laplace transform is given by 
\[
l_{\mu _{\mathtt{G}}^{\sigma }}(\varphi )=\exp \left( -\left\langle \log
(1-\varphi )\right\rangle _{\sigma }\right) ,\;1>\varphi \in \mathcal{D}. 
\]
This measure can be seen as a special case of compound Poisson measure $\mu
_{\mathtt{CP}}$ for a specific choice of the measure $\rho $ used in the
definition of $\mu _{\mathtt{CP}}$, see Section \ref{2eq66} - (\ref{2eq67})
for details. From this point of view, of course, all structure may be
implemented on Gamma space. The question that still remains is to find
intrinsic expressions for all these operators on Gamma space as found in
Poisson space.

The most intriguing feature of Gamma space we found is its Fock type
structure. As in the Poisson case it is possible to choose a transformation $%
\alpha $ on $\mathcal{D}$ such that the normalized exponentials $e_{\mu _{%
\mathtt{G}}^{\sigma }}^{\sigma }(\varphi ;\cdot )$ produce a complete system
of orthogonal polynomials, the so called system of generalized Laguerre
polynomials. It leads to a Fock type realization of Gamma space as 
\[
L^{2}(\mu _{\mathtt{G}}^{\sigma })\simeq \bigoplus_{n=0}^{\infty }\mathrm{Exp%
}_{n}^{\mathtt{G}}L^{2}(\sigma )=\mathrm{Exp}^{\mathtt{G}}L^{2}(\sigma ),
\]
where $\mathrm{Exp}_{{n}}^{\mathtt{G}}L^{2}(\sigma )\subset \mathrm{Exp}%
_{n}L^{2}(\sigma )$ is a quasi-$n$-particle subspace of $\mathrm{Exp}^{%
\mathtt{G}}L^{2}(\sigma ).$ The point here is that the scalar product in $%
\mathrm{Exp}_{n}^{\mathtt{G}}L^{2}(\sigma )$ turns out to be different of
the standard one given by $L^{2}(\sigma )^{\widehat{\otimes }n}.$ As a
result the space $\mathrm{Exp}^{\mathtt{G}}L^{2}(\sigma )$ has a novel $n$%
-particle structure which is essentially different from traditional Fock
picture.

\section{Poisson analysis\label{ref21}}

Throughout this section we consider the measure space $(\QTR{mathbb}{R}^{d},%
\mathcal{B}(\QTR{mathbb}{R}^{d}),\sigma )$ where $\sigma $ is a non-atomic $%
\sigma $-finite measure. Below we denote by $\mathcal{D}:=\mathcal{D}(%
\QTR{mathbb}{R}^{d})$ and by $\mathcal{D}^{\prime }:=\mathcal{D}^{\prime }(%
\QTR{mathbb}{R}^{d})$ the classical Schwartz spaces and by $\mathcal{M}(%
\QTR{mathbb}{R}^{d})\subset \mathcal{D}^{\prime }$ the set of all positive
Radon measures on $($the set of all positive Radon measures on $(%
\QTR{mathbb}{R}^{d},\mathcal{B}(\QTR{mathbb}{R}^{d})).$

\subsection{The configuration space over $\QTR{mathbb}{R}^{d}$ \label{2eq69}}

The configuration space $\Gamma _{\QTR{mathbb}{R}^{d}}=:\Gamma $ over $%
\QTR{mathbb}{R}^{d}$ is defined as the set of all locally finite subsets
(configurations) in $\QTR{mathbb}{R}^{d},$ i.e., 
\begin{equation}
\Gamma :=\left\{ \gamma \subset \QTR{mathbb}{R}^{d}\,|\,\left| \gamma \cap
K\right| <\infty \;\mathrm{for\;any\;compact\;}K\subset \QTR{mathbb}{R}%
^{d}\right\} .  \label{2eq70}
\end{equation}
Here $\left| A\right| $ denotes the cardinality of a set $A.$

We can identify any $\gamma \in \Gamma $ with the corresponding sum of Dirac
measures, namely 
\begin{equation}
\Gamma \ni \gamma \rightarrow \sum_{x\in \gamma }\varepsilon _{x}\left(
dy\right) =:d\gamma \left( y\right) \in \mathcal{M}_{p}(\QTR{mathbb}{R}%
^{d})\subset \mathcal{M}(\QTR{mathbb}{R}^{d}),  \label{2eq23}
\end{equation}
where $\mathcal{M}_{p}(\QTR{mathbb}{R}^{d})$ denotes the set of all positive
integer valued measures $($or Radon point measures$)$ over $\mathcal{B}(%
\QTR{mathbb}{R}^{d}). $

The space $\Gamma $ can be endowed with the relative topology as a closed
subset of the space $\mathcal{M}_{p}(\QTR{mathbb}{R}^{d})$ on $\mathcal{B}(%
\QTR{mathbb}{R}^{d})$ with the vague topology, i.e., a sequence of measures $%
(\mu _{n})_{n\in \QTR{mathbb}{N}}$ converge in the vague topology to $\mu $
if and only if for any $f\in C_{0}(\QTR{mathbb}{R}^{d})$ (the set of all
continuous functions with compact support) we have 
\[
\int_{\QTR{mathbb}{R}^{d}}f\left( x\right) d\mu _{n}\left( x\right) 
\stackunder{n\rightarrow \infty }{\longrightarrow }\int_{\QTR{mathbb}{R}%
^{d}}f\left( x\right) d\mu \left( x\right) . 
\]
Then for any $f\in C_{0}(\QTR{mathbb}{R}^{d})$ we have a continuous
functional 
\[
\Gamma \ni \gamma \mapsto \left\langle \gamma ,f\right\rangle :=\left\langle
f\right\rangle _{\gamma }=\int_{\QTR{mathbb}{R}^{d}}f\left( x\right) d\gamma
\left( x\right) =\sum_{x\in \gamma }f\left( x\right) \in \QTR{mathbb}{R}. 
\]
Conversely, such functionals generate the topology of the space $\Gamma .$

Hence we have the following chain 
\[
\Gamma \subset \mathcal{M}(\QTR{mathbb}{R}^{d})\subset \mathcal{D}^{\prime
}:=\mathcal{D}^{\prime }\left( \QTR{mathbb}{R}^{d}\right) . 
\]
The Borel $\sigma $-algebra on $\Gamma ,$ $\mathcal{B}(\Gamma ),$ is
generated by sets of the form 
\begin{equation}
C_{\Lambda ,n}=\left\{ \gamma \in \Gamma \,|\,\left| \gamma \cap \Lambda
\right| =n\right\} ,  \label{2eq57}
\end{equation}
where $\Lambda \in \mathcal{B}(\QTR{mathbb}{R}^{d})$ bounded$,$ $n\in 
\QTR{mathbb}{N}$, see e.g. \cite{GGV75} and for any $\Lambda \in \mathcal{B}(%
\QTR{mathbb}{R}^{d})$ and all $n\in \QTR{mathbb}{N}$ the set $C_{\Lambda ,n}$
is a Borel set of $\Gamma .$ Sets of the form (\ref{2eq57}) are called
cylinder sets.

For any $B\subset \QTR{mathbb}{R}^{d}$ we introduce a function $N_{B}:\Gamma
\rightarrow \QTR{mathbb}{N}$ such that 
\begin{equation}
N(B)(\gamma )=\left| \gamma \cap B\right| ,\,\gamma \in \Gamma .
\label{2eq79}
\end{equation}
Then $\mathcal{B}(\Gamma )$ is the minimal $\sigma $-algebra with which all
the functions $\{N_{B}\,|\,B\in \mathcal{B}(\QTR{mathbb}{R}^{d})$ bounded$\}$
are measurable.

\subsection{The Poisson measure and its properties}

The Poisson measure $\pi _{\sigma }$ (with intensity measure $\sigma $) on ($%
\Gamma ,\mathcal{B}(\Gamma )$) may be defined in different ways, here we
give two convenient characterizations of $\pi _{\sigma }.$

\begin{definition}[Laplace transform]
The Laplace transform of $\pi _{\sigma }$ is given by 
\begin{equation}
l_{\pi _{\sigma }}\left( \varphi \right) =\int_{\Gamma }\exp \left(
\left\langle \gamma ,\varphi \right\rangle \right) d\pi _{\sigma }\left(
\gamma \right) =\exp \left( \int_{\QTR{mathbb}{R}^{d}}\left( e^{\varphi
\left( x\right) }-1\right) d\sigma \left( x\right) \right) ,  \label{2eq1}
\end{equation}
where $\varphi \in \mathcal{D}$, see e.g. \cite{KMM78} and 
\cite[Chap.~III Sec.~4]{GV68}.
\end{definition}

\begin{remark}
The right hand side of (\ref{2eq1}) defines, via Minlos' theorem, the
measure $\pi _{\sigma }$ on $(\mathcal{D}^{\prime },\mathcal{C}_{\sigma }(%
\mathcal{D}^{\prime })),$ but an additional analysis shows that the support
of the measure $\pi _{\sigma }$ is $\Gamma \subset \mathcal{D}^{\prime },$
see e.g. \cite{Ka74}, \cite{Ka83} and \cite{KMM78}, hence $\pi _{\sigma }$
can be considered as a measure on $\Gamma .$
\end{remark}

Let $f:\QTR{mathbb}{R}^{d}\times \Gamma \rightarrow \QTR{mathbb}{R}$ be such
that $f\geq 0$ and measurable with respect to $\mathcal{B}(\QTR{mathbb}{R}%
^{d})\times \mathcal{B}(\Gamma )$. Define 
\begin{eqnarray*}
F\left( \gamma \right) &\mbox{$:=$}&\,\left\langle \gamma ,f\left( \cdot
,\gamma \right) \right\rangle \\
&=&\int_{\QTR{mathbb}{R}^{d}}f\left( x,\gamma \right) d\gamma \left( x\right)
\\
&=&\sum_{x\in \gamma }f\left( x,\gamma \right) .
\end{eqnarray*}
Then $\pi _{\sigma }$ is characterized by 
\begin{eqnarray}
\int_{\Gamma }F\left( \gamma \right) d\pi _{\sigma }\left( \gamma \right) &%
\mbox{$:=$}&\,\int_{\Gamma }\int_{\QTR{mathbb}{R}^{d}}f\left( x,\gamma
\right) d\gamma \left( x\right) d\pi _{\sigma }\left( \gamma \right) 
\nonumber \\
&=&\,\,\int_{\QTR{mathbb}{R}^{d}}\int_{\Gamma }f\left( x,\gamma +\varepsilon
_{x}\right) d\pi _{\sigma }\left( \gamma \right) d\sigma \left( x\right) .
\label{2eq3}
\end{eqnarray}
Equality (\ref{2eq3}) is known as Mecke identity, see e.g. \cite{Me67} and 
\cite{NZ79}.

\subsection{The Fock space isomorphism of Poisson space\label{2eq61}}

Let us consider the following transformation on $\mathcal{D},$ $\alpha :%
\mathcal{D}\rightarrow \mathcal{D}$ defined by 
\[
\alpha \left( \varphi \right) \left( x\right) =\log \left( 1+\varphi \left(
x\right) \right) ,\;-1<\varphi \in \mathcal{D},\,x\in \QTR{mathbb}{R}^{d}. 
\]
As easily can be seen $\alpha \left( 0\right) =0$ and $\alpha $ is
holomorphic in some neighborhood $\mathcal{U}_{\alpha }$ of zero. Using this
transformation we introduce the normalized exponential $e_{\pi _{\sigma
}}^{\alpha }(\cdot ,\cdot )$ which is holomorphic on a neighborhood of zero $%
\mathcal{U}_{\alpha }^{\prime }\subset \mathcal{U}_{\alpha }\subset \mathcal{%
D}.$ For $\varphi \in \mathcal{U}_{\alpha }^{\prime },$ $\gamma \in \Gamma $
we set 
\begin{eqnarray}
e_{\pi _{\sigma }}^{\alpha }\left( \varphi ,\gamma \right) &\mbox{$:=$}&\,%
\frac{\exp \left( \left\langle \gamma ,\alpha \left( \varphi \right)
\right\rangle \right) }{l_{\pi _{\sigma }}\left( \alpha \left( \varphi
\right) \right) }  \nonumber \\
&=&\,\exp \left( \left\langle \gamma ,\log \left( 1+\varphi \right)
\right\rangle -\left\langle \varphi \right\rangle _{\sigma }\right) ,
\label{2eq2}
\end{eqnarray}
where $\left\langle \varphi \right\rangle _{\sigma }:=\int_{\QTR{mathbb}{R}%
^{d}}\varphi \left( x\right) d\sigma \left( x\right) .$

We use the holomorphy of $e_{\pi _{\sigma }}^{\alpha }(\cdot ,\gamma )$ on a
neighborhood of zero to expand it in power series which, with Cauchy's
inequality, polarization identity and kernel theorem, give us the following
result 
\begin{equation}
e_{\pi _{\sigma }}^{\alpha }\left( \varphi ,\gamma \right)
=\sum_{n=0}^{\infty }\dfrac{1}{n!}\left\langle P_{n}^{\pi _{\sigma },\alpha
}\left( \gamma \right) ,\varphi ^{\otimes n}\right\rangle ,\,\varphi \in 
\mathcal{U}_{\alpha }^{\prime }\subset \mathcal{U}_{\alpha },\,\gamma \in
\Gamma ,  \label{2eq59}
\end{equation}
where $P_{n}^{\pi _{\sigma },\alpha }:\Gamma \rightarrow \mathcal{D}^{\prime 
\widehat{\otimes }n}.$ $\{P_{n}^{\pi _{\sigma },\alpha }\left( \cdot \right)
=:C_{n}^{\sigma }\left( \cdot \right) \,|\,n\in \QTR{mathbb}{N}_{0}\}$ is
called the system of \textbf{generalized Charlier kernels }on Poisson space $%
(\Gamma ,\mathcal{B}(\Gamma ),\pi _{\sigma }).$ From (\ref{2eq59}) it
follows immediately that for any $\varphi ^{(n)}\in \mathcal{D}^{\widehat{%
\otimes }n},$ $n\in \QTR{mathbb}{N}_{0}$ the function 
\[
\Gamma \ni \gamma \mapsto \left\langle C_{n}^{\sigma }\left( \gamma \right)
,\varphi ^{(n)}\right\rangle 
\]
is a polynomial of the order $n$ on $\Gamma .$ The system of functions 
\[
\left\{ C_{n}^{\sigma }\left( \varphi ^{(n)}\right) \left( \gamma \right)
:=\left\langle C_{n}^{\sigma }\left( \gamma \right) ,\varphi ^{\left(
n\right) }\right\rangle ,\,\forall \varphi ^{(n)}\in \mathcal{D}^{\widehat{%
\otimes }n},\,n\in \QTR{mathbb}{N}_{0}\right\} 
\]
is called the system of \textbf{generalized Charlier polynomials} for the
Poisson measure $\pi _{\sigma }.$

\begin{proposition}
For any $\varphi ^{\left( n\right) }\in \mathcal{D}^{\widehat{\otimes }n}$
and $\psi ^{\left( m\right) }\in \mathcal{D}^{\widehat{\otimes }m}$ we have 
\[
\int_{\Gamma }\left\langle C_{n}^{\sigma }\left( \gamma \right) ,\varphi
^{\left( n\right) }\right\rangle \left\langle C_{m}^{\sigma }\left( \gamma
\right) ,\psi ^{\left( m\right) }\right\rangle d\pi _{\sigma }\left( \gamma
\right) =\delta _{nm}n!\left( \varphi ^{\left( n\right) },\psi ^{\left(
n\right) }\right) _{L^{2}\left( \sigma ^{\otimes n}\right) }. 
\]
\end{proposition}

\noindent \textbf{Proof.\ }Let $\varphi ^{\left( n\right) },\psi ^{\left(
m\right) }$ be given as in the proposition and such that $\varphi ^{\left(
n\right) }=\varphi ^{\otimes n},$ $\psi ^{\left( m\right) }=\psi ^{\otimes
m} $. Then for $z_{1},z_{2}\in \QTR{mathbb}{C,}$ and taking into account (%
\ref{2eq1}) and (\ref{2eq2}) we have 
\begin{eqnarray}
&&\int_{\Gamma }e_{\pi _{\sigma }}^{\alpha }(z_{1}\varphi ,\gamma )e_{\pi
_{\sigma }}^{\alpha }(z_{2}\psi ,\gamma )d\pi _{\sigma }\left( \gamma \right)
\nonumber \\
&=&\exp \left( -\left\langle z_{1}\varphi +z_{2}\psi \right\rangle _{\sigma
}\right) \int_{\Gamma }\exp \left( \left\langle \gamma ,\log \left( \left(
1+z_{1}\varphi \right) \left( 1+z_{2}\psi \right) \right) \right\rangle
\right) d\pi _{\sigma }\left( \gamma \right)  \nonumber \\
&=&\exp \left( -\left\langle z_{1}\varphi +z_{2}\psi \right\rangle _{\sigma
}\right)  \nonumber \\
&&\cdot \exp \left( \int_{\QTR{mathbb}{R}^{d}}\left( \exp \left( \log \left(
\left( 1+z_{1}\varphi \right) \left( 1+z_{2}\psi \right) \right) \right)
-1\right) d\sigma \right)  \nonumber \\
&=&\exp \left( z_{1}z_{2}\left( \varphi ,\psi \right) _{L^{2}(\sigma
)}\right)  \nonumber \\
&=&\sum_{n=0}^{\infty }\frac{1}{n!}z_{1}^{n}z_{2}^{n}\left( \varphi
^{\otimes n},\psi ^{\otimes n}\right) _{L^{2}(\sigma ^{\otimes n})}.
\label{2eq15}
\end{eqnarray}
On the other hand 
\begin{eqnarray}
&&\int_{\Gamma }e_{\pi _{\sigma }}^{\alpha }(z_{1}\varphi ,\gamma )e_{\pi
_{\sigma }}^{\alpha }(z_{2}\psi ,\gamma )d\pi _{\sigma }\left( \gamma \right)
\nonumber \\
&=&\sum_{n,m=0}^{\infty }\frac{z_{1}^{n}z_{2}^{m}}{n!m!}\int_{\Gamma
}\left\langle C_{n}^{\sigma }\left( \gamma \right) ,\varphi ^{\otimes
n}\right\rangle \left\langle C_{m}^{\sigma }\left( \gamma \right) ,\psi
^{\otimes m}\right\rangle d\pi _{\sigma }\left( \gamma \right) .
\label{2eq16}
\end{eqnarray}
Then a comparison of coefficients between (\ref{2eq15}) and (\ref{2eq16}),
with polarization identity and linearity, gives the above result.\hfill $%
\blacksquare $

\begin{remark}
This proposition gives us the possibility to extend - in the $L^{2}(\pi
_{\sigma })$ sense - the class of $\langle C_{n}^{\sigma }\left( \gamma
\right) ,\varphi ^{\otimes n}\rangle $-functions to include kernels from the
so-called $n$-particle Fock space over $L^{2}\left( \sigma \right) .$
\end{remark}

We define the Fock space as the Hilbert sum 
\[
\mathrm{Exp}L^{2}\left( \sigma \right) :=\bigoplus_{n=0}^{\infty }\mathrm{Exp%
}_{n}L^{2}\left( \sigma \right) 
\]
where $\mathrm{Exp}_{n}L^{2}\left( \sigma \right) :=L^{2}\left( \sigma
\right) _{\QTR{mathbb}{C}}^{\widehat{\otimes }n}$ and we put by definition $%
\mathrm{Exp}_{0}L^{2}\left( \sigma \right) :=\QTR{mathbb}{C}$.

For any $F\in L^{2}\left( \pi _{\sigma }\right) $ there exists a sequence $%
(f^{\left( n\right) })_{n=0}^{\infty }\in \mathrm{Exp}L^{2}\left( \sigma
\right) $ such that 
\begin{equation}
F\left( \gamma \right) =\sum_{n=0}^{\infty }\left\langle C_{n}^{\sigma
}\left( \gamma \right) ,f^{\left( n\right) }\right\rangle  \label{2eq47}
\end{equation}
and moreover 
\begin{equation}
\left\| F\right\| _{L^{2}\left( \pi _{\sigma }\right)
}^{2}=\sum_{n=0}^{\infty }n!\left| f^{\left( n\right) }\right| _{L^{2}\left(
\sigma ^{\otimes n}\right) }^{2},  \label{2eq46}
\end{equation}
where the r.h.s. of (\ref{2eq46}) coincides with the square of the norm in $%
\mathrm{Exp}L^{2}\left( \sigma \right) .$ And vice versa, any series of the
form (\ref{2eq47}) with coefficients $(f^{\left( n\right) })_{n=0}^{\infty
}\in \mathrm{Exp}L^{2}\left( \sigma \right) $ gives a function from $%
L^{2}\left( \pi _{\sigma }\right) .$ As a result we have the well-known
isomorphism $I_{\sigma }$ between $L^{2}\left( \pi _{\sigma }\right) $ and $%
\mathrm{Exp}L^{2}\left( \sigma \right) .$

Now we introduce the action of the annihilation and creation operators in
the Fock space $\mathrm{Exp}L^{2}(\sigma )$, see e.g. \cite[Appendix A.2]
{HKPS93} and \cite{RS75}. Consider $f^{\left( n\right) }\in \mathrm{Exp}%
_{n}L^{2}\left( \sigma \right) $ of the form 
\begin{equation}
f^{\left( n\right) }=\widehat{\otimes }_{i=1}^{n}f_{i},\;f_{i}\in
L^{2}\left( \sigma \right) ,\,i=1,\ldots ,n.  \label{2eq42}
\end{equation}
Then the action of the annihilation operator $a^{-}(\varphi ),$ $\varphi \in 
\mathcal{D},$ on $f^{\left( n\right) }$ is defined as follows: 
\[
a^{-}\left( \varphi \right) f^{\left( n\right) }:=\sum_{j=1}^{n}\left\langle
\varphi ,f_{j}\right\rangle \widehat{\otimes }_{{\QATOP{i=1}{i\neq j}}%
}^{n}f_{i}\in \mathrm{Exp}_{n-1}L^{2}\left( \sigma \right) .
\]
This definition is independent of the particular representation of $%
f^{\left( n\right) }$ in (\ref{2eq42}), hence $a^{-}(\varphi )f^{\left(
n\right) }$ is well-defined. Moreover this definition can be extended by
linearity to a dense subspace of $\mathrm{Exp}_{n}L^{2}\left( \sigma \right) 
$ consisting of finite linear combinations of elements of the form (\ref
{2eq42}). One easily finds the following inequality for such elements 
\begin{equation}
\left| a^{-}(\varphi )f^{\left( n\right) }\right| \leq \sqrt{n}\left|
\varphi \right| \left| f^{\left( n\right) }\right| ,  \label{2eq48}
\end{equation}
which shows that the extension of $a^{-}(\varphi )$ to $\mathrm{Exp}%
_{n}L^{2}\left( \sigma \right) $ as a bounded operator exists. Consider the
dense subspace $\mathrm{Exp}_{0}L^{2}\left( \sigma \right) $ of $\mathrm{Exp}%
L^{2}\left( \sigma \right) $ consisting of those sequences $\{f^{\left(
n\right) },n\in \QTR{mathbb}{N}_{0}\}$ which only have a finite numbers of
non-vanishing entries. The bound (\ref{2eq48}) allow us to extend $%
a^{-}(\varphi ),$ $\varphi \in \mathcal{D}$, component-wise to $\mathrm{Exp}%
_{0}L^{2}\left( \sigma \right) $ which, therefore, give us a densely defined
operator on $\mathrm{Exp}L^{2}\left( \sigma \right) $ denoted again by $%
a^{-}(\varphi ).$ So the adjoint operator of $a^{-}(\varphi )$ exists, which
we denote by $a^{+}(\varphi ),$ and call creation operator. The action of
the creation operator on elements $f^{\left( n\right) }\in $ $\mathrm{Exp}%
_{n}L^{2}\left( \sigma \right) $ is given by 
\[
a^{+}(\varphi )f^{\left( n\right) }=\varphi \widehat{\otimes }f^{\left(
n\right) }\in \mathrm{Exp}_{n+1}L^{2}\left( \sigma \right) .
\]
For the creation operator we also have an estimate 
\[
\left| a^{+}(\varphi )f^{\left( n\right) }\right| \leq \sqrt{n+1}\left|
\varphi \right| \left| f^{\left( n\right) }\right| .
\]
As before, this estimate give us the possibility to deduce that in the same $%
a^{+}(\varphi )$ is densely defined on $\mathrm{Exp}L^{2}\left( \sigma
\right) .$

For latter use we introduce a vector $\mathrm{Exp}\psi ,$ $\psi \in
L^{2}(\sigma )$ as 
\[
\mathrm{Exp}\psi =\left( \frac{1}{n!}\psi ^{\otimes n}\right) _{n=0}^{\infty
} 
\]
which is called the coherent state corresponding to the one-particle state $%
\psi $ (or exponential vector corresponding to $\psi \in L^{2}(\sigma )$).
For any set $\mathcal{L}\subset L^{2}(\sigma )$ which is total in $%
L^{2}(\sigma )$ the set of coherent states $\{\mathrm{Exp}\psi \,|\,\psi \in 
\mathcal{L}\}\subset \mathrm{Exp}L^{2}(\sigma )$ is also total in $\mathrm{%
Exp}L^{2}(\sigma )$, see e.g. \cite[Chap. 2]{G72} and \cite{BK88}. We note
that $e_{\pi _{\sigma }}^{\alpha }(\psi ,\cdot )$ is nothing as the coherent
state in the Fock space picture, for any $\psi \in \mathcal{D},$ $\psi >-1,$
we have 
\[
L^{2}\left( \pi _{\sigma }\right) \ni e_{\pi _{\sigma }}^{\alpha }(\psi
,\cdot )=\sum_{n=0}^{\infty }\frac{1}{n!}\left\langle C_{n}^{\sigma }\left(
\cdot \right) ,\psi ^{\otimes n}\right\rangle \mapsto \mathrm{Exp}\psi \in 
\mathrm{Exp}L^{2}\left( \sigma \right) . 
\]

The action of the annihilation operator $a^{-}(\varphi )$ on $\mathrm{Exp}%
\psi $ is given by 
\[
a^{-}\left( \varphi \right) \mathrm{Exp}\psi =\left( \varphi ,\psi \right)
_{L^{2}\left( \sigma \right) }\mathrm{Exp}\psi . 
\]

\subsection{Annihilation operator on Poisson space}

Let us introduce a set of smooth cylinder functions $\mathcal{F}%
C_{b}^{\infty }\left( \mathcal{D},\Gamma \right) $ $($dense in $L^{2}(\pi
_{\sigma }))$ which consists of all functions of the form 
\[
f\left( \gamma \right) =F\left( \left\langle \gamma ,\varphi
_{1}\right\rangle ,\ldots ,\left\langle \gamma ,\varphi _{N}\right\rangle
\right) ,\;\gamma \in \Gamma , 
\]
where the generating directions $\varphi _{1},\ldots ,\varphi _{N}\in 
\mathcal{D},$ and $F$ (generating function for $f$) is from $C_{b}^{\infty }(%
\QTR{mathbb}{R}^{N})$ ($C^{\infty }$-functions on $\QTR{mathbb}{R}^{d}$ with
bounded derivatives).

\begin{definition}
We define the \textbf{\emph{Poissonian gradient }}$\nabla ^{\mathtt{P}}$ as
a mapping 
\[
\nabla ^{\mathtt{P}}:\mathcal{F}C_{b}^{\infty }\left( \mathcal{D},\Gamma
\right) \longrightarrow L^{2}\left( \pi _{\sigma }\right) \otimes
L^{2}\left( \sigma \right) 
\]
given by 
\[
\left( \nabla ^{\mathtt{P}}f\right) \left( \gamma ,x\right) =f\left( \gamma
+\varepsilon _{x}\right) -f\left( \gamma \right) ,\;\gamma \in \Gamma
,\;x\in \QTR{mathbb}{R}^{d}.
\]
\end{definition}

Let us mention that the operation 
\[
\Gamma \ni \gamma \mapsto \gamma +\varepsilon _{x}\in \Gamma 
\]
is well-defined because of the property: 
\[
\pi _{\sigma }\left\{ \gamma \in \Gamma \,|\,x\in \gamma \right\}
=0,\;\forall x\in \QTR{mathbb}{R}^{d}. 
\]
The fact that $\mathcal{F}C_{b}^{\infty }\left( \mathcal{D},\Gamma \right)
\ni f\mapsto \nabla ^{\mathtt{P}}f\in L^{2}\left( \pi _{\sigma }\right)
\otimes L^{2}\left( \sigma \right) $ arises from the use of the Hilbert
space $L^{2}\left( \sigma \right) $ as a tangent space at any point $\gamma
\in \Gamma .$

\begin{remark}
To produce differential structure we need linear structure where the measure
have support $\Gamma .\,$If we consider $\pi _{\sigma }$ on $\mathcal{D}%
^{\prime }$ then $\pi _{\sigma }(\xi +\varphi )\perp \pi _{\sigma }(\xi ),$
see e.g. \cite{GGV75}, therefore integration by parts and adjoint of
operators are not available. This is the reason why we embed $\Gamma $ in $%
\mathcal{D}^{\prime }.$
\end{remark}

\begin{proposition}
\label{2eq50}For any non-negative function $h\in \mathrm{Dom}((\nabla ^{%
\mathtt{P}})^{*})$ such that $h\in L^{1}(\pi _{\sigma })\otimes L^{1}(\sigma
)$ the following equality holds 
\begin{equation}
\left( \left( \nabla ^{\mathtt{P}}\right) ^{*}h\right) \left( \gamma \right)
=\int_{\QTR{mathbb}{R}^{d}}h\left( \gamma -\varepsilon _{x},x\right) d\gamma
\left( x\right) -\int_{\QTR{mathbb}{R}^{d}}h\left( \gamma ,x\right) d\sigma
\left( x\right) .  \label{2eq5}
\end{equation}
\end{proposition}

\noindent \textbf{Proof.\ }Let $f\in \mathrm{Dom}(\nabla ^{\mathtt{P}})$ be
given. Then we use the Mecke identity (\ref{2eq3}) to compute $\left( \nabla
^{\mathtt{P}}f,h\right) _{L^{2}\left( \pi _{\sigma }\right) \otimes
L^{2}\left( \sigma \right) }$ as follows: 
\begin{eqnarray*}
&&\left( \nabla ^{\mathtt{P}}f,h\right) _{L^{2}\left( \pi _{\sigma }\right)
\otimes L^{2}\left( \sigma \right) } \\
&=&\int_{\QTR{mathbb}{R}^{d}}\int_{\Gamma }\left( f\left( \gamma
+\varepsilon _{x}\right) -f\left( \gamma \right) \right) h\left( \gamma
,x\right) d\pi _{\sigma }\left( \gamma \right) d\sigma \left( x\right) \\
&=&\int_{\QTR{mathbb}{R}^{d}}\int_{\Gamma }f\left( \gamma +\varepsilon
_{x}\right) h\left( \gamma ,x\right) d\pi _{\sigma }\left( \gamma \right)
d\sigma \left( x\right) \\
&&-\int_{\QTR{mathbb}{R}^{d}}\int_{\Gamma }f\left( \gamma \right) h\left(
\gamma ,x\right) d\pi _{\sigma }\left( \gamma \right) d\sigma \left( x\right)
\\
&=&\int_{\Gamma }f\left( \gamma \right) \stackunder{\left( \left( \nabla ^{%
\mathtt{P}}\right) ^{*}h\right) \left( \gamma \right) }{\underbrace{\left[
\int_{\QTR{mathbb}{R}^{d}}h\left( \gamma -\varepsilon _{x},x\right) d\gamma
\left( x\right) -\int_{\QTR{mathbb}{R}^{d}}h\left( \gamma ,x\right) d\sigma
\left( x\right) \right] }}d\pi _{\sigma }\left( \gamma \right) .
\end{eqnarray*}
\hfill $\blacksquare $

Now we are going to give an internal description of the annihilation
operator.

The directional derivative is then defined as 
\begin{eqnarray}
\left( \nabla _{\varphi }^{\mathtt{P}}f\right) \left( \gamma \right)
&=&\left( \left( \nabla ^{\mathtt{P}}f\right) \left( \gamma ,\cdot \right)
,\varphi \left( \cdot \right) \right) _{L^{2}\left( \sigma \right) } 
\nonumber \\
&=&\int_{\QTR{mathbb}{R}^{d}}\left( f\left( \gamma +\varepsilon _{x}\right)
-f\left( \gamma \right) \right) \varphi \left( x\right) d\sigma \left(
x\right)  \label{2eq6}
\end{eqnarray}
for any $\varphi \in \mathcal{D}.$ Of course the operator 
\[
\nabla _{\varphi }^{\mathtt{P}}:\mathcal{F}C_{b}^{\infty }\left( \mathcal{D}%
,\Gamma \right) \longrightarrow L^{2}\left( \pi _{\sigma }\right) 
\]
is closable in $L^{2}\left( \pi _{\sigma }\right) .$

\begin{proposition}
\label{2eq43}The closure of $\nabla _{\varphi }^{\mathtt{P}}$ coincides with
the image under $I_{\sigma }$ of the annihilation operator $a^{-}(\varphi )$
in $\mathrm{Exp}L^{2}\left( \sigma \right) ,$ i.e., $I_{\sigma
}a^{-}(\varphi )I_{\sigma }^{-1}=\nabla _{\varphi }^{\mathtt{P}}.$
\end{proposition}

\noindent \textbf{Proof.\ }To prove this proposition it is enough to show
this equality of operators in a total set in the core of the annihilation
operator. Let $\psi \in \mathcal{U}_{\alpha }^{\prime }$ be given, then
having in mind (\ref{2eq6}) and (\ref{2eq2}) it follows that 
\begin{eqnarray}
\left( \nabla _{\varphi }^{\mathtt{P}}e_{\pi _{\sigma }}^{\alpha }(\psi
;\cdot )\right) \left( \gamma \right) &=&\int_{\QTR{mathbb}{R}^{d}}\left(
e_{\pi _{\sigma }}^{\alpha }(\psi ;\gamma +\varepsilon _{x})-e_{\pi _{\sigma
}}^{\alpha }(\psi ;\gamma )\right) \varphi \left( x\right) d\sigma \left(
x\right)  \nonumber \\
&=&e_{\pi _{\sigma }}^{\alpha }\left( \psi ;\gamma \right) \int_{%
\QTR{mathbb}{R}^{d}}\left( \exp \left( \left\langle \varepsilon _{x},\log
\left( 1+\psi \right) \right\rangle \right) -1\right) \varphi \left(
x\right) d\sigma \left( x\right)  \nonumber \\
&=&\left( \psi ,\varphi \right) _{L^{2}\left( \sigma \right) }e_{\pi
_{\sigma }}^{\alpha }\left( \psi ;\gamma \right) .  \label{2eq44}
\end{eqnarray}
On the other hand since $I_{\sigma }^{-1}e_{\pi _{\sigma }}^{\alpha }(\psi
;\gamma )=\mathrm{Exp}\psi $ it follows that 
\[
a^{-}(\varphi )\mathrm{Exp}\psi =\left( \varphi ,\psi \right) _{L^{2}\left(
\sigma \right) }\mathrm{Exp}\psi . 
\]
Hence if we apply $I_{\sigma }$ to this vector we just obtain the same
result as (\ref{2eq44}) which had to be proven.\hfill $\blacksquare $

\subsection{Creation operator on Poisson space\label{2eq63}}

\begin{proposition}
For any $\varphi \in \mathcal{D},$ $g\in \mathrm{Dom}(I_{\sigma
}a^{+}(\varphi )I_{\sigma }^{-1}),$ where $a^{+}(\varphi )\,$is the creation
operator in $\mathrm{Exp}L^{2}(\sigma ),$ the following equality holds 
\begin{eqnarray}
\left( \left( \nabla _{\varphi }^{\mathtt{P}}\right) ^{*}g\right) \left(
\gamma \right) &=&\int_{\QTR{mathbb}{R}^{d}}g\left( \gamma -\varepsilon
_{x}\right) \varphi \left( x\right) d\gamma \left( x\right) -g\left( \gamma
\right) \int_{\QTR{mathbb}{R}^{d}}\varphi \left( x\right) d\sigma \left(
x\right)  \nonumber \\
&=&\left( g\left( \gamma -\varepsilon _{\cdot }\right) ,\varphi \left( \cdot
\right) \right) _{L^{2}(\gamma )}-g\left( \gamma \right) \left\langle
\varphi \right\rangle _{\sigma }.  \label{2eq4}
\end{eqnarray}
\end{proposition}

\begin{remark}
In terms of chaos decomposition of $g\in \mathrm{Dom}((\nabla _{\varphi }^{%
\mathtt{P}})^{*})$ the equality (\ref{2eq4}) was established in \cite{NV95}.
We give an independent proof of (\ref{2eq4}), which is based on the results
on absolute continuity of Poisson measures, see e.g. \cite{Sk57} and \cite
{T90}.
\end{remark}

\noindent \textbf{Proof.\ 1.} First we give a version of the proof of (\ref
{2eq4}) which uses the Mecke identity.

It follows from (\ref{2eq6}) that 
\begin{eqnarray}
\left( \nabla _{\varphi }^{\mathtt{P}}f,g\right) _{L^{2}\left( \pi _{\sigma
}\right) } &=&\int_{\Gamma }\left( \left( \nabla ^{\mathtt{P}}f\right)
\left( \gamma ,\cdot \right) ,\varphi \left( \cdot \right) \right)
_{L^{2}\left( \sigma \right) }g\left( \gamma \right) d\pi _{\sigma }\left(
\gamma \right)  \nonumber \\
&=&\left( \left( \nabla ^{\mathtt{P}}f\right) \left( \cdot ,\cdot \right)
,g\left( \cdot \right) \varphi \left( \cdot \right) \right) _{L^{2}\left(
\pi _{\sigma }\right) \otimes L^{2}\left( \sigma \right) }.  \label{2eq39}
\end{eqnarray}
Whence using Proposition \ref{2eq50} we obtain 
\begin{eqnarray*}
\left( \left( \nabla _{\varphi }^{\mathtt{P}}\right) ^{*}g\right) \left(
\gamma \right) &=&\left( \left( \nabla ^{\mathtt{P}}\right) ^{*}g\varphi
\right) \left( \gamma \right) \\
&=&\int_{\QTR{mathbb}{R}^{d}}g\left( \gamma -\varepsilon _{x}\right) \varphi
\left( x\right) d\gamma \left( x\right) -g\left( \gamma \right) \int_{%
\QTR{mathbb}{R}^{d}}\varphi \left( x\right) d\sigma \left( x\right) \\
&=&\left( g\left( \gamma -\varepsilon _{\cdot }\right) ,\varphi \left( \cdot
\right) \right) _{L^{2}(\gamma )}-g\left( \gamma \right) \left\langle
\varphi \right\rangle _{\sigma },
\end{eqnarray*}
which proves (\ref{2eq4}).

\noindent \textbf{2.} Alternatively we give an independent prove of (\ref
{2eq4}) based on absolute continuity of Poisson measure.

Let $\eta \in \mathcal{D}$ be such that $\eta \left( x\right) >-1,\,\forall
x\in \QTR{mathbb}{R}^{d}.$ Denote by $\sigma _{\eta }$ the measure on $%
\QTR{mathbb}{R}^{d}$ having density with respect to $\sigma ,$%
\begin{equation}
\frac{d\sigma _{\eta }}{d\sigma }\left( x\right) =1+\eta \left( x\right) .
\label{2eq9}
\end{equation}

\begin{lemma}
The Poisson measure $\pi _{\sigma }$ and $\pi _{\sigma _{\eta }}$ on $%
(\Gamma ,\mathcal{B}(\Gamma ))$ are mutually absolutely continuous and the
Radon-Nikodym derivative $\frac{d\pi _{\sigma _{\eta }}}{d\pi _{\sigma }}%
(\gamma )$ coincides with the normalized exponential, i.e., 
\[
\frac{d\pi _{\sigma _{\eta }}}{d\pi _{\sigma }}(\gamma )=e_{\pi _{\sigma
}}^{\alpha }\left( \eta ;\gamma \right) =\exp \left( \left\langle \gamma
,\log \left( 1+\eta \right) \right\rangle -\left\langle \eta \right\rangle
_{\sigma }\right) . 
\]
\end{lemma}

\noindent \textbf{Proof.\ }Let $\eta \in \mathcal{D}$ be such that $\eta
\left( x\right) >-1,\,\forall x\in \QTR{mathbb}{R}^{d}.$ Then the Laplace
transform of $\pi _{\sigma _{\eta }},$ given by (\ref{2eq1}), 
\begin{eqnarray*}
&&\int_{\Gamma }\exp \left( \left\langle \gamma ,\varphi \right\rangle
\right) d\pi _{\sigma _{\eta }}\left( \gamma \right) \\
&=&\exp \left( \int_{\QTR{mathbb}{R}^{d}}\left( e^{\varphi \left( x\right)
}-1\right) \left( 1+\eta \left( x\right) \right) d\sigma \left( x\right)
\right) \\
&=&e^{-\left\langle \eta \right\rangle _{\sigma }}\exp \left( \int_{%
\QTR{mathbb}{R}^{d}}\left( e^{\varphi \left( x\right) +\log \left( 1+\eta
\left( x\right) \right) }-1\right) d\sigma \left( x\right) \right) \\
&=&\int_{\Gamma }\exp \left( \left\langle \gamma ,\varphi \right\rangle
\right) \stackunder{d\pi _{\sigma _{\eta }}\left( \gamma \right) }{%
\underbrace{\exp \left( \left\langle \gamma ,\log \left( 1+\eta \right)
\right\rangle -\left\langle \eta \right\rangle _{\sigma }\right) d\pi
_{\sigma }\left( \gamma \right) }}.
\end{eqnarray*}
\hfill $\blacksquare $

\hspace{0in}In order to proof (\ref{2eq4}) it suffices to verify the
equality 
\begin{equation}
\left( \nabla _{\varphi }^{\mathtt{P}}f,g\right) _{L^{2}\left( \pi _{\sigma
}\right) }=\int_{\Gamma }f\left( \gamma \right) \left[ \left( g\left( \gamma
-\varepsilon _{\cdot }\right) ,\varphi \left( \cdot \right) \right)
_{L^{2}(\gamma )}-g\left( \gamma \right) \left\langle \varphi \right\rangle
_{\sigma }\right] d\pi _{\sigma }\left( \gamma \right)  \label{2eq7}
\end{equation}
for $f(\gamma )=e_{\pi _{\sigma }}^{\alpha }(\psi ;\gamma ),\;g(\gamma
)=e_{\pi _{\sigma }}^{\alpha }(\eta ;\gamma ),\;\psi ,\eta $ belong to a
neighborhood of zero $\mathcal{U}\subset \mathcal{D},$ because the coherent
states $\mathrm{Exp}\psi ,\;\psi \in \mathcal{U}$ span a common core for the
annihilation and creation operators.

\begin{lemma}
For any $\varphi \in \mathcal{D}$ and for all $\psi ,\eta $ in a
neighborhood of zero $\mathcal{U}_{\alpha }^{\prime }\subset \mathcal{D},$
the following equality holds 
\begin{equation}
\left( \nabla _{\varphi }^{\mathtt{P}}e_{\pi _{\sigma }}^{\alpha }(\psi
;\cdot ),e_{\pi _{\sigma }}^{\alpha }(\eta ;\cdot )\right) _{L^{2}\left( \pi
_{\sigma }\right) }=\left( \psi ,\varphi \right) _{L^{2}\left( \sigma
\right) }\exp \left( \left( \psi ,\eta \right) _{L^{2}\left( \sigma \right)
}\right) .  \label{2eq45}
\end{equation}
\end{lemma}

\noindent \textbf{Proof.\ }Taking in account (\ref{2eq44}) we compute the
right hand side of (\ref{2eq45}) to be 
\begin{eqnarray}
&&\left( \nabla _{\varphi }^{\mathtt{P}}e_{\pi _{\sigma }}^{\alpha }(\psi
;\cdot ),e_{\pi _{\sigma }}^{\alpha }(\eta ;\cdot )\right) _{L^{2}\left( \pi
_{\sigma }\right) }  \nonumber \\
&=&\left( \psi ,\varphi \right) _{L^{2}\left( \sigma \right) }\exp \left(
-\left\langle \psi +\eta \right\rangle _{\sigma }\right) \int_{\Gamma }\exp
\left( \left\langle \gamma ,\log \left( \left( 1+\psi \right) \left( 1+\eta
\right) \right) \right\rangle \right) d\pi _{\sigma }\left( \gamma \right) 
\nonumber \\
&=&\left( \psi ,\varphi \right) _{L^{2}\left( \sigma \right) }\exp \left(
-\left\langle \psi +\eta \right\rangle _{\sigma }\right)  \nonumber \\
&&\cdot \exp \left( \int_{\QTR{mathbb}{R}^{d}}\left( \psi \left( x\right)
+\eta \left( x\right) +\psi \left( x\right) \eta \left( x\right) \right)
d\sigma \left( x\right) \right)  \nonumber \\
&=&\left( \psi ,\varphi \right) _{L^{2}\left( \sigma \right) }\exp \left(
\left( \psi ,\eta \right) _{L^{2}\left( \sigma \right) }\right) ,
\label{2eq11}
\end{eqnarray}
which proves the statement of the lemma.\hfill $\blacksquare $

Further, the r.h.s. of (\ref{2eq7}) can be rewritten as follows 
\begin{eqnarray}
&&\int_{\Gamma }e_{\pi _{\sigma }}^{\alpha }\left( \psi ;\gamma \right)
\left[ \sum_{x\in \gamma }\exp \left( \left\langle \gamma -\varepsilon
_{x},\log \left( 1+\eta \right) \right\rangle -\left\langle \eta
\right\rangle _{\sigma }\right) \varphi \left( x\right) \right.  \nonumber \\
&&\left. \left. -e_{\pi _{\sigma }}^{\alpha }\left( \eta ;\gamma \right)
\left\langle \varphi \right\rangle _{\sigma }\right) \QATOP{\QATOP {} {}}{%
\QATOP {} {}}\!\!\right] d\pi _{\sigma }\left( \gamma \right)  \nonumber \\
&=&\int_{\Gamma }e_{\pi _{\sigma }}^{\alpha }\left( \psi ;\gamma \right)
e_{\pi _{\sigma }}^{\alpha }\left( \eta ;\gamma \right) \left\langle \tfrac{%
\varphi }{1+\eta }\right\rangle _{\gamma }d\pi _{\sigma }\left( \gamma
\right)  \nonumber \\
&&-\left\langle \varphi \right\rangle _{\sigma }\exp \left( \left( \psi
,\eta \right) _{L^{2}(\sigma )}\right) .  \label{2eq8}
\end{eqnarray}
Let us state the following useful lemma.

\begin{lemma}
\begin{enumerate}
\item  $\left\langle \psi \right\rangle _{\sigma _{\eta }}=\left\langle \psi
\right\rangle _{\sigma }+\left( \psi ,\eta \right) _{L^{2}(\sigma
)},\;\forall \psi ,\eta \in \mathcal{D}.$

\item  $e_{\pi _{\sigma _{\eta }}}^{\alpha }(\psi ;\gamma )=\exp \left(
-\left( \psi ,\eta \right) _{L^{2}(\sigma )}\right) e_{\pi _{\sigma
}}^{\alpha }(\psi ;\gamma ),\;\forall \psi \in \mathcal{U}\subset \mathcal{D}%
.$

\item  \label{2eq10}$\left\langle \gamma ,\frac{\psi }{1+\eta }\right\rangle
=\left\langle C_{1}^{\sigma _{\eta }}\left( \gamma \right) ,\frac{\psi }{%
1+\eta }\right\rangle +\left\langle \frac{\psi }{1+\eta }\right\rangle
_{\sigma _{\eta }}.$
\end{enumerate}
\end{lemma}

\noindent \textbf{Proof.\ }The non-trivial step is \ref{2eq10}. Let us
denote for simplicity $\frac{\psi }{1+\eta }=:\xi $%
\begin{eqnarray*}
\left\langle C_{1}^{\sigma _{\eta }}\left( \gamma \right) ,\xi \right\rangle
&=&\left. \frac{d}{dt}e_{\pi _{\sigma _{\eta }}}^{\alpha }\left( t\xi
;\gamma \right) \right| _{t=0} \\
&=&\left. \frac{d}{dt}\exp \left( \left\langle \gamma ,\log \left( 1+t\xi
\right) \right\rangle -\left\langle t\xi \right\rangle _{\sigma _{\eta
}}\right) \right| _{t=0} \\
&=&\left. \frac{d}{dt}\sum_{x\in \gamma }\log \left( 1+t\xi \left( x\right)
\right) -\left\langle t\xi \right\rangle _{\sigma _{\eta }}\right| _{t=0} \\
&=&\left\langle \gamma ,\xi \right\rangle -\left\langle \xi \right\rangle
_{\sigma _{\eta }}.
\end{eqnarray*}
\hfill $\blacksquare $

\noindent Now the rest of the proof follows from the previous lemma and (\ref
{2eq8}), i.e., 
\begin{eqnarray*}
&&\int_{\Gamma }e_{\pi _{\sigma }}^{\alpha }\left( \psi ;\gamma \right)
e_{\pi _{\sigma }}^{\alpha }\left( \eta ;\gamma \right) \left\langle \tfrac{%
\varphi }{1+\eta }\right\rangle _{\gamma }d\pi _{\sigma }\left( \gamma
\right) -\left\langle \varphi \right\rangle _{\sigma }\exp \left( \left(
\psi ,\eta \right) _{L^{2}(\sigma )}\right) \\
&=&\exp \left( \left( \psi ,\eta \right) _{L^{2}(\sigma )}\right) \left(
\int_{\Gamma }e_{\pi _{\sigma }}^{\alpha }\left( \psi ;\gamma \right)
\left\langle \tfrac{\varphi }{1+\eta }\right\rangle _{\gamma }d\pi _{\sigma
_{\eta }}\left( \gamma \right) -\left\langle \varphi \right\rangle _{\sigma
}\right) \\
&=&\exp \left( \left( \psi ,\eta \right) _{L^{2}(\sigma )}\right) \left(
\int_{\Gamma }e_{\pi _{\sigma }}^{\alpha }\left( \psi ;\gamma \right)
\left\langle C_{1}^{\sigma _{\eta }}\left( \gamma \right) ,\tfrac{\varphi }{%
1+\eta }\right\rangle d\pi _{\sigma _{\eta }}\left( \gamma \right) \right. \\
&&\left. +\left\langle \tfrac{\varphi }{1+\eta }\right\rangle _{\sigma
_{\eta }}\int_{\Gamma }e_{\pi _{\sigma _{\eta }}}^{\alpha }\left( \psi
;\gamma \right) d\pi _{\sigma _{\eta }}\left( \gamma \right) -\left\langle
\varphi \right\rangle _{\sigma }\right) \\
&=&\exp \left( \left( \psi ,\eta \right) _{L^{2}(\sigma )}\right) \left(
\psi ,\frac{\varphi }{1+\eta }\right) _{L^{2}(\sigma _{\eta })} \\
&=&\left( \psi ,\varphi \right) _{L^{2}(\sigma )}\exp \left( \left( \psi
,\eta \right) _{L^{2}(\sigma )}\right) ,
\end{eqnarray*}
which is the same as (\ref{2eq11}). This completes the proof.\hfill $%
\blacksquare $

\begin{remark}
\label{2eq62}The operator $(\nabla _{\varphi }^{\mathtt{P}})^{*}\,$plays the
role of creation operator since \linebreak $a^{+}(\varphi )^{n}1=\varphi
^{\otimes n},$ i.e., 
\begin{equation}
\left( \left( \nabla _{\varphi }^{\mathtt{P}}\right) ^{*n}1\right) \left(
\gamma \right) =\left\langle C_{n}^{\sigma }\left( \gamma \right) ,\varphi
^{\otimes n}\right\rangle .  \label{2eq56}
\end{equation}
\end{remark}

\section{Compound Poisson measures}

\subsection{Definition and properties}

This section is devoted to study the compound Poisson measures $\mu _{%
\mathtt{CP}}$ on $(\mathcal{D}^{\prime },\mathcal{B}(\mathcal{D}^{\prime
})). $ Firstly we recall the L\'{e}vy canonical representation of all
possible generalized white noise measures $\mu $ on $(\mathcal{D}^{\prime },%
\mathcal{B}(\mathcal{D}^{\prime }))$, see \cite{GV68}, \cite{H70}, and \cite
{AW95}. These measures are defined by the characteristic functional of the
form 
\begin{equation}
C_{\mu }(\varphi )=\exp \left[ ia(\varphi ,1)-\frac{b^{2}|\varphi |^{2}}{2}%
+\int_{\QTR{mathbb}{R}}(e^{is\varphi }-1-\frac{is\varphi }{1+s^{2}},1)d\beta
(s)\right] ,  \label{2eq38}
\end{equation}
where $\varphi \in \mathcal{D}$, $a,b\in \QTR{mathbb}{R}$, the measure $%
\beta $ is such that $\beta (\{0\})=0$, $\int_{\QTR{mathbb}{R}%
}s^{2}/(1+s^{2})d\beta (s)<\infty $, $(\cdot ,\cdot )$ and $|\cdot |$ denote
the inner product and the norm in $L^{2}(\QTR{mathbb}{R}^{d})$,
respectively. We take into account that such a measure is in general the
convolution of a Gaussian and non Gaussian measures. We will be interested
in the non Gaussian part of this class, i.e., in (\ref{2eq38}) $b=0$.
Furthermore, assume that 
\[
\int_{\QTR{mathbb}{R}}s^{2}d\beta (s)<\infty ,\;\int_{-1}^{1}|s|d\beta
(s)<\infty . 
\]
Then one can use the Kolmogorov canonical representation of $C_{\mu
}(\varphi )$, as in e.g.~\cite{GV68} 
\[
C_{\mu }(\varphi )=\exp \left[ ia(\varphi ,1)+\int_{\QTR{mathbb}{R}%
}(e^{is\varphi }-1,1)d\beta (s)\right] ,\;\varphi \in \mathcal{D}. 
\]

Let us define the compound Poisson measures on $(\mathcal{D}^{\prime },%
\mathcal{B}(\mathcal{D}^{\prime }))$. Let $\rho \,$be a measure on $(%
\QTR{mathbb}{R},\mathcal{B}(\QTR{mathbb}{R}))$ (finite or $\sigma $-finite)
having all moments finite and such that $\rho (\{0\})=0$. In addition let $%
\sigma $ be a non-atomic $\sigma $-finite measure on $(\QTR{mathbb}{R}^{d},%
\mathcal{B}(\QTR{mathbb}{R}^{d}))$. We denote 
\[
\psi _{\rho }\left( u\right) :=\int_{\QTR{mathbb}{R}}\left( e^{su}-1\right)
d\rho \left( s\right) ,\;s\in \QTR{mathbb}{R}. 
\]

\begin{definition}
A measure $\mu _{\mathtt{CP}}$ on $(\mathcal{D}^{\prime },\mathcal{B}(%
\mathcal{D}^{\prime }))$ is called a \textbf{\emph{compound Poisson measure}}
with Kolmogorov characteristic $\psi _{\rho }$ if its Laplace transform is
given by, as e.g. \cite{GGV75} 
\begin{eqnarray}
l_{\mu _{\mathtt{CP}}}\left( \varphi \right) &=&\int_{\mathcal{D}^{\prime
}}\exp \left( \left\langle \omega ,\varphi \right\rangle \right) d\mu _{%
\mathtt{CP}}\left( \omega \right)  \nonumber \\
&=&\exp \left( \int_{\QTR{mathbb}{R}^{d}}\psi _{\rho }\left( \varphi \left(
x\right) \right) d\sigma \left( x\right) \right)  \nonumber \\
&=&\exp \left( \int_{\QTR{mathbb}{R}^{d}}\int_{\QTR{mathbb}{R}}\left(
e^{s\varphi \left( x\right) }-1\right) d\rho \left( s\right) d\sigma \left(
x\right) \right) ,\,\varphi \in \mathcal{D}.  \label{2eq22}
\end{eqnarray}
\end{definition}

\begin{proposition}
\begin{enumerate}
\item  \label{2eq30}Assume that $\rho $ satisfies the analyticity property: 
\begin{equation}
\exists C>0:\forall n\in \QTR{mathbb}{N}\;\int_{\QTR{mathbb}{R}}\left|
s\right| ^{n}d\rho \left( s\right) <C^{n}n!.  \label{2eq72}
\end{equation}
Then the Laplace transform of $\mu _{\mathtt{CP}}$ is holomorphic at $0\in 
\mathcal{D}_{\QTR{mathbb}{C}}$.

\item  \label{2eq31}Let $\rho (\QTR{mathbb}{R})<\infty .$ Then 
\[
\mu _{\mathtt{CP}}\left( \Omega \right) :=\mu _{\mathtt{CP}}\left( \left\{
\sum_{x\in \gamma }s_{x}\varepsilon _{x}\in \mathcal{D}^{\prime }|s_{x}\in 
\mathrm{supp\,}\rho ,\,\gamma \in \Gamma \right\} \right) =1, 
\]
where $\Gamma =\Gamma _{\QTR{mathbb}{R}^{d}}$ is the same as defined in
Subsection \ref{2eq69}-(\ref{2eq70}).

\item  \label{2eq32}Let $\rho (\QTR{mathbb}{R})=\infty .$ Then 
\begin{equation}
\mu _{\mathtt{CP}}\left( \Omega \right) :=\mu _{\mathtt{CP}}\left( \left\{
\sum_{x\in \gamma _{c}}s_{x}\varepsilon _{x}\in \mathcal{D}^{\prime
}|s_{x}\in \mathrm{supp\,}\rho ,\,\gamma _{c}\in \Gamma _{c}\right\} \right)
=1,  \label{2eq34}
\end{equation}
where $\Gamma _{c}$ is the collection of all locally countable subets in $%
\QTR{mathbb}{R}^{d}$.
\end{enumerate}
\end{proposition}

\textbf{Proof.\ }\ref{2eq30}. By (\ref{2eq72}) the Kolmogorov characteristic 
$\psi _{\rho }$ is holomorphic on some neighborhood of $0\in \QTR{mathbb}{C}%
\QTR{mathbb}{.}$ Then by (\ref{2eq22}) the Laplace transform $l_{\mu _{%
\mathtt{CP}}}$ of $\mu _{\mathtt{CP}}$ is holomorphic in some neighborhood
of zero $\mathcal{U}\subset \mathcal{D}_{\QTR{mathbb}{C}}.$

\noindent \ref{2eq31}, \ref{2eq32}. At first assume $d=1.$ Then $\mu _{%
\mathtt{CP}}$ corresponds to the distributional derivative of the compound
Poisson process $\xi _{t}$ and statements \ref{2eq31}, \ref{2eq32} follow
immediately from the properties of the paths of this process. Namely, almost
every path $\xi _{t}$ of compound Poisson process is right continuous step
function with the jumps from $\mathrm{supp}\rho .$ If $\rho $ is finite a
measure, then any finite interval contains only finite number of the points
of discontinuities of $\xi _{t}$ (in this case $\xi $ is called a
generalized Poisson process). For infinite measure $\rho $ the set of
discontinuities of $\xi _{t}$ is locally countable, see e.g. \cite{Ta67} and 
\cite{K93}.

For $d>1$ the statements \ref{2eq31}, \ref{2eq32} follow from the analogous
results of the theory of Poisson measures, see e.g. \cite{Ka74}, \cite{Ka83}
and \cite{KMM78}.\hfill $\blacksquare $

\begin{remark}
Assume that $\rho $ is a probability measure on $(\QTR{mathbb}{R,}\mathcal{B}%
(\QTR{mathbb}{R}))$. Let $\{\xi _{k},k\geq 1\}$ be a sequence of independent
identically $\rho $-distributed random variables and $N=\{N_{t},t\geq 0\}$
be the standard Poisson process independent of $\{\xi _{k},k\geq 1\}$. Then $%
\mu _{\mathtt{CP}}$ is generated by the distributional derivative of the
compound Poisson process 
\[
\xi _{t}=\sum_{k=1}^{N_{t}}\xi _{k}. 
\]
Notice that we don't consider in this paper measures $\mu $ corresponding to
the distributional derivatives of doubly stochastic Poisson processes and
fields.
\end{remark}

\subsection{The isomorphism between Poisson and compound Poisson spaces\label%
{2eq64}}

Let us define the measure $\widehat{\sigma }$ on $(\QTR{mathbb}{R}^{d+1},%
\mathcal{B}(\QTR{mathbb}{R}^{d+1}))$ as the product of the measures $\rho $
and $\sigma ,$ i.e., 
\[
d\widehat{\sigma }\left( \widehat{x}\right) :=d\rho \left( s\right) d\sigma
\left( x\right) ,\;\widehat{x}=\left( s,x\right) \in \QTR{mathbb}{R}\times 
\QTR{mathbb}{R}^{d} 
\]
and consider on $(\Gamma _{\QTR{mathbb}{R}^{d+1}},\mathcal{B}(\Gamma _{%
\QTR{mathbb}{R}^{d+1}}))$ the Poisson measure $\pi _{\widehat{\sigma }}$
with intensity measure $\widehat{\sigma }$. According to (\ref{2eq1}) this
means that the Laplace transform of $\pi _{\widehat{\sigma }}$ has the form 
\begin{eqnarray}
l_{\pi _{\widehat{\sigma }}}(\widehat{\varphi }) &=&\int_{\Gamma _{%
\QTR{mathbb}{R}^{d+1}}}\exp (\langle \widehat{\gamma },\widehat{\varphi }%
\rangle )d\pi _{\widehat{\sigma }}(\widehat{\gamma })  \nonumber \\
&=&\exp \left( \int_{\QTR{mathbb}{R}^{d+1}}(e^{\widehat{\varphi }(\widehat{x}%
)}-1)d\widehat{\sigma }(\widehat{x})\right) ,\;\widehat{\varphi }\in 
\mathcal{D}(\QTR{mathbb}{R}^{d+1}).  \label{2eq33}
\end{eqnarray}

The intensity measure $\widehat{\sigma }$ has the following property: for
any $x\in \QTR{mathbb}{R}^{d}$, $\Delta \in \mathcal{B}(\QTR{mathbb}{R})$
such that $\rho (\Delta )<\infty $%
\[
\widehat{\sigma }(\Delta \times \{x\})=\rho (\Delta )\sigma (\{x\})=0. 
\]
This property yields that $\pi _{\widehat{\sigma }}$ is concentrated on a
smaller set than $\Gamma _{\QTR{mathbb}{R}^{d+1}}$. Namely, let us define $%
\widehat{\Gamma }\subset \Gamma _{\QTR{mathbb}{R}^{d+1}}$ as follows 
\begin{equation}
\widehat{\Gamma }:=\left\{ \widehat{\gamma }\in \Gamma _{\QTR{mathbb}{R}%
^{d+1}}|\widehat{\gamma }=\sum_{\widehat{x}_{i}\in \widehat{\gamma }%
}\varepsilon _{\widehat{x}_{i}},\;\widehat{x}_{i}=\left( s_{i},x_{i}\right)
\in \QTR{mathbb}{R}\times \QTR{mathbb}{R}^{d},\,x_{i}\neq x_{j},\;i\neq
j\right\} .  \label{2eq73}
\end{equation}

\begin{proposition}
\label{2eq81}The measure $\pi _{\widehat{\sigma }}$ is concentrated on the
set $\widehat{\Gamma }\in \mathcal{B}(\Gamma _{\QTR{mathbb}{R}^{d+1}})$.
\end{proposition}

\noindent \textbf{Proof.} One can deduce this result from the theory of
Point processes, see e.g.~\cite[Chap.~1]{KMM78}, \cite{Ka83}, and \cite{K93}.

\begin{remark}
Let $d=1$. The measure $\pi _{\widehat{\sigma }}$ corresponds to the
distributional derivative of the independently marked Poisson process with
the intensity measure $\sigma $ and space of marks $(\QTR{mathbb}{R,\rho })$
(for more details on marked processes see, e.g.~\cite{K93} and \cite{BL95}).
For $d>1$ there exists analogous connection between $\pi _{\widehat{\sigma }}
$ and independently marked Poisson fields with the same intensity and
marking.
\end{remark}

It follows from (\ref{2eq38}) that the Laplace transform $l_{\pi _{\widehat{%
\sigma }}}$ is well defined for $\widehat{\varphi }\left( s,x\right)
=p\left( s\right) \varphi \left( x\right) $ where $p\left( s\right)
=\sum_{k=0}^{m}p_{k}s^{k}$ ($p_{0}\neq 0$ for finite $\rho $ and $p_{0}=0$
for infinite $\rho $) is a polynomial and $\varphi \in \mathcal{D}$ (cf. 
\cite{LRS95}). Let us put $\widehat{\varphi }\left( s,x\right) =s\varphi
\left( x\right) ,$ $\varphi \in \mathcal{D}$ in (\ref{2eq33}). Then by (\ref
{2eq22}) we obtain 
\[
l_{\mu _{\mathtt{CP}}}\left( \varphi \right) =l_{\pi _{\widehat{\sigma }%
}}\left( s\varphi \right) ,\;\varphi \in \mathcal{D}. 
\]

Then using (\ref{2eq34}) it follows that the compound Poisson measure $\mu _{%
\mathtt{CP}}$ is the image of $\pi _{\widehat{\sigma }}$ under the
transformation $\Sigma :\widehat{\Gamma }\rightarrow \Sigma \widehat{\Gamma }%
=\Omega \subset \mathcal{D}^{\prime }$ given by 
\begin{equation}
\widehat{\Gamma }\ni \widehat{\gamma }\mapsto \left( \Sigma \widehat{\gamma }%
\right) \left( \cdot \right) =\Sigma \left( \sum_{\widehat{x}_{i}\in 
\widehat{\gamma }}\varepsilon _{\widehat{x}_{i}}\right) \left( \cdot \right)
:=\sum_{\left( s_{i},x_{i}\right) \in \widehat{\gamma }}s_{i}\varepsilon
_{x_{i}}\left( \cdot \right) \in \Omega \subset \mathcal{D}^{\prime },
\label{2eq35}
\end{equation}
i.e., $\forall B\in \mathcal{B}(\mathcal{D}^{\prime })$ 
\[
\mu _{\mathtt{CP}}\left( B\right) =\mu _{\mathtt{CP}}\left( B\cap \Omega
\right) =\pi _{\widehat{\sigma }}\left( \Sigma ^{-1}\left( B\cap \Omega
\right) \right) , 
\]
where $\Sigma ^{-1}\Delta $ is the pre-image of the set $\Delta .$

The latter equality may be rewritten in the following form 
\[
\int_{\mathcal{D}^{\prime }}1\!\!1_{B}\left( \omega \right) d\mu _{\mathtt{CP%
}}\left( \omega \right) =\int_{\Omega }1\!\!1_{B}\left( \omega \right) d\mu
_{\mathtt{CP}}\left( \omega \right) =\int_{\widehat{\Gamma }%
}1\!\!1_{B}\left( \Sigma \widehat{\gamma }\right) d\pi _{\widehat{\sigma }%
}\left( \widehat{\gamma }\right) , 
\]
which is nothing than the well known change of variable formula for the
Lebesgue integral. Namely, for any $h\in L^{1}(\Omega ,\mu _{\mathtt{CP}})$
the function $h\circ \Sigma \in L^{1}(\widehat{\Gamma },\pi _{\widehat{%
\sigma }})$ and 
\begin{equation}
\int_{\Omega }h\left( \omega \right) d\mu _{\mathtt{CP}}\left( \omega
\right) =\int_{\widehat{\Gamma }}h\left( \Sigma \widehat{\gamma }\right)
d\pi _{\widehat{\sigma }}\left( \widehat{\gamma }\right) .  \label{2eq36}
\end{equation}

\begin{remark}
It is worth noting that there exists on $\Omega $ an inverse map $\Sigma
^{-1}:\Omega \rightarrow \widehat{\Gamma }.$ And we obtain that $\pi _{%
\widehat{\sigma }}$ on $\widehat{\Gamma }$ is the image of $\mu _{\mathtt{CP}%
}$ on $\Omega $ under the map $\Sigma ^{-1},$ i.e., $\forall \widehat{C}\in 
\mathcal{B}(\widehat{\Gamma }),$ $\pi _{\widehat{\sigma }}(\widehat{C})=\mu
_{\mathtt{CP}}(\Sigma \widehat{C})$ or after rewriting 
\[
\int_{\widehat{\Gamma }}1\!\!1_{\widehat{C}}\left( \widehat{\gamma }\right)
d\pi _{\widehat{\sigma }}\left( \widehat{\gamma }\right) =\int_{\Omega
}1\!\!1_{\Sigma \widehat{C}}\left( \omega \right) d\mu _{\mathtt{CP}}\left(
\omega \right) =\int_{\Omega }1\!\!1_{\widehat{C}}\left( \Sigma ^{-1}\omega
\right) d\mu _{\mathtt{CP}}\left( \omega \right) . 
\]
As before we easily can write the corresponding change of variables formula,
namely $\forall \widehat{f}\in L^{1}(\widehat{\Gamma },\pi _{\widehat{\sigma 
}})$ the function $\widehat{f}\circ \Sigma ^{-1}\in L^{1}(\Omega ,\mu _{%
\mathtt{CP}})$ and 
\[
\int_{\widehat{\Gamma }}\widehat{f}\left( \widehat{\gamma }\right) d\pi _{%
\widehat{\sigma }}\left( \widehat{\gamma }\right) =\int_{\Omega }\widehat{f}%
\left( \Sigma ^{-1}\omega \right) d\mu _{\mathtt{CP}}\left( \omega \right) . 
\]
\end{remark}

So we construct a unitary isomorphism $U_{\Sigma }$ between the Poisson
space $L^{2}(\pi _{\widehat{\sigma }})=L^{2}(\widehat{\Gamma },\pi _{%
\widehat{\sigma }})$ and the compound Poisson space $L^{2}(\mu _{\mathtt{CP}%
})=L^{2}(\Omega ,\mu _{\mathtt{CP}}).$ Namely, 
\[
L^{2}\left( \Omega ,\mu _{\mathtt{CP}}\right) \ni h\mapsto U_{\Sigma
}h:=h\circ \Sigma \in L^{2}\left( \widehat{\Gamma },\pi _{\widehat{\sigma }%
}\right) 
\]
and 
\begin{equation}
L^{2}\left( \widehat{\Gamma },\pi _{\widehat{\sigma }}\right) \ni \widehat{f}%
\mapsto U_{\Sigma }^{-1}\widehat{f}=\widehat{f}\circ \Sigma ^{-1}\in
L^{2}\left( \Omega ,\mu _{\mathtt{CP}}\right) .  \label{2eq71}
\end{equation}
The isometry of $U_{\Sigma }$ and $U_{\Sigma }^{-1}$ follows from (\ref
{2eq36}).

As a result we have established the following proposition.

\begin{proposition}
The map $U_{\Sigma }$ is a unitary isomorphism between the Poisson space and
the compound Poisson space.
\end{proposition}

\begin{remark}
In the space $L^{2}(\pi _{\widehat{\sigma }})$ we have a basis of
generalized Charlier polynomials, annihilation and creation operators etc.
Now we can use the unitary isomorphism $U_{\Sigma }$ in order to transport
the Fock structure from $L^{2}(\pi _{\widehat{\sigma }})$ to $L^{2}(\mu _{%
\mathtt{CP}}).$
\end{remark}

\subsection{Annihilation and creation operators on compound Poisson space%
\label{2eq65}}

Let $\bigtriangledown _{\widehat{\varphi }}^{\mathtt{P}},$ $%
(\bigtriangledown _{\widehat{\varphi }}^{\mathtt{P}})^{*},$ $\widehat{%
\varphi }\in \mathcal{D}(\QTR{mathbb}{R}^{d+1})$ be the annihilation and
creation operators on Poisson space $L^{2}(\pi _{\widehat{\sigma }}).$ Their
images under $U_{\Sigma }$%
\begin{equation}
U_{\Sigma }^{-1}\bigtriangledown _{\widehat{\varphi }}^{\mathtt{P}}U_{\Sigma
},\;\;U_{\Sigma }^{-1}\left( \bigtriangledown _{\widehat{\varphi }}^{\mathtt{%
P}}\right) ^{*}U_{\Sigma }  \label{2eq51}
\end{equation}
play the role of annihilation and creation operators in compound Poisson
space $L^{2}(\mu _{\mathtt{CP}}).$ Let us calculate the actions of (\ref
{2eq51}).

The set of smooth cylinder functions $\mathcal{F}C_{b}^{\infty }(\mathcal{D}%
,\Omega ),$ $($dense in $L^{2}(\mu _{\mathtt{CP}}))$ consists of all
functions of the form 
\begin{eqnarray}
h\left( \omega \right) &=&H\left( \left\langle \omega ,\varphi
_{1}\right\rangle ,\ldots ,\left\langle \omega ,\varphi _{N}\right\rangle
\right)  \nonumber \\
&=&H\left( \left\langle \Sigma ^{-1}\omega ,s\varphi _{1}\right\rangle
,\ldots ,\left\langle \Sigma ^{-1}\omega ,s\varphi _{N}\right\rangle \right)
,  \label{2eq41}
\end{eqnarray}
where (generating directions) $\varphi _{1},\ldots ,\varphi _{N}\in \mathcal{%
D}$ and $H$ (generating function for $h$) is from $C_{b}^{\infty }(%
\QTR{mathbb}{R}^{N})$. Whence it follows that 
\[
\mathcal{F}C_{b}^{\infty }\left( \mathcal{D},\Omega \right) =U_{\Sigma }^{-1}%
\mathcal{F}C_{b}^{\infty }\left( \mathcal{D}\left( \QTR{mathbb}{R}%
^{d+1}\right) ,\widehat{\Gamma }\right) . 
\]

By (\ref{2eq6}) for any $\widehat{f}\in \mathcal{F}C_{b}^{\infty }(\mathcal{D%
}(\QTR{mathbb}{R}^{d+1}),\widehat{\Gamma })$ we have 
\begin{equation}
\left( \nabla _{\widehat{\varphi }}^{\mathtt{P}}\widehat{f}\right) \left( 
\widehat{\gamma }\right) =\int_{\QTR{mathbb}{R}^{d+1}}\left( \widehat{f}%
\left( \widehat{\gamma }+\varepsilon _{\widehat{x}}\right) -\widehat{f}%
\left( \widehat{\gamma }\right) \right) \widehat{\varphi }\left( \widehat{x}%
\right) d\widehat{\sigma }\left( \widehat{x}\right) .  \label{2eq55}
\end{equation}

\begin{proposition}
For any $h\in \mathcal{F}C_{b}^{\infty }(\mathcal{D},\Omega )$ the operator $%
U_{\Sigma }^{-1}\bigtriangledown _{\widehat{\varphi }}^{\mathtt{P}}U_{\Sigma
}$ has the following form 
\[
\left( U_{\Sigma }^{-1}\bigtriangledown _{\widehat{\varphi }}^{\mathtt{P}%
}U_{\Sigma }h\right) \left( \omega \right) =\int_{\QTR{mathbb}{R}^{d}}\int_{%
\QTR{mathbb}{R}}\left( h\left( \omega +s\varepsilon _{x}\right) -h\left(
\omega \right) \right) \widehat{\varphi }\left( s,x\right) d\rho \left(
s\right) d\sigma \left( x\right) . 
\]
\end{proposition}

\noindent \textbf{Proof.\ }Let $h\in \mathcal{F}C_{b}^{\infty }(\mathcal{D}%
,\Omega )$ be given and denote $U_{\Sigma }h=h\circ \Sigma =:\widehat{h}$
and $\Sigma ^{-1}w=:\widehat{\gamma }$. Taking into account (\ref{2eq51})
and (\ref{2eq41}) we obtain 
\begin{eqnarray}
\left( U_{\Sigma }^{-1}\bigtriangledown _{\widehat{\varphi }}^{\mathtt{P}%
}U_{\Sigma }h\right) \left( \omega \right) &=&\left( \bigtriangledown _{%
\widehat{\varphi }}^{\mathtt{P}}\widehat{h}\right) \left( \widehat{\gamma }%
\right)  \nonumber \\
&=&\int_{\QTR{mathbb}{R}^{d+1}}\left( \widehat{h}\left( \widehat{\gamma }%
+\varepsilon _{\widehat{x}}\right) -\widehat{h}\left( \widehat{\gamma }%
\right) \right) \widehat{\varphi }\left( \widehat{x}\right) d\widehat{\sigma 
}\left( \widehat{x}\right) .  \label{2eq52}
\end{eqnarray}
Now we use the definition of $\widehat{h},$ the additivity of the map $%
\Sigma $ and the obvious equality $\Sigma \varepsilon _{\widehat{x}%
}=s\varepsilon _{x}$ for $\widehat{x}=(s,x)$; with this (\ref{2eq52}) turns
out to be 
\begin{eqnarray}
&&\int_{\QTR{mathbb}{R}^{d+1}}\left( h\left( \Sigma \left( \widehat{\gamma }%
+\varepsilon _{\widehat{x}}\right) \right) -h\left( \Sigma \widehat{\gamma }%
\right) \widehat{\varphi }\left( \widehat{x}\right) \right) d\widehat{\sigma 
}\left( \widehat{x}\right)  \nonumber \\
&=&\int_{\QTR{mathbb}{R}^{d+1}}\left( h\left( \omega +s\varepsilon
_{x}\right) -h\left( \omega \right) \right) \widehat{\varphi }\left( 
\widehat{x}\right) d\widehat{\sigma }\left( \widehat{x}\right) .
\label{2eq53}
\end{eqnarray}
The result of the proposition follows then by definition of $\widehat{\sigma 
}.$\hfill $\blacksquare $

Putting $\widehat{\varphi }\left( \widehat{x}\right) =\phi \left( s\right)
\varphi \left( x\right) $ in (\ref{2eq53}) we obtain 
\begin{eqnarray*}
&&\left( U_{\Sigma }^{-1}\bigtriangledown _{\phi \varphi }^{\mathtt{P}%
}U_{\Sigma }h\right) \left( \omega \right) \\
&=&\int_{\QTR{mathbb}{R}^{d}}\left( \int_{\QTR{mathbb}{R}}\left( h\left(
\omega +s\varepsilon _{x}\right) -h\left( \omega \right) \right) \phi \left(
s\right) d\rho \left( s\right) \right) \varphi \left( x\right) d\sigma
\left( x\right) .
\end{eqnarray*}
Let us note that by (\ref{2eq38}) we can admit not only bounded functions $%
\phi \left( s\right) $ but also polynomials. For finite $\rho $ and $\phi
\equiv 1$ we have the following formula for the \textbf{annihilation
operator }$\bigtriangledown _{\varphi }^{\mathtt{CP}}$ in compound Poisson
space $L^{2}(\Omega ,\mu _{\mathtt{CP}}):$%
\begin{eqnarray}
\left( \bigtriangledown _{\varphi }^{\mathtt{CP}}h\right) \left( \omega
\right) &\mbox{$:=$}&\left( U_{\Sigma }^{-1}\bigtriangledown _{\varphi }^{%
\mathtt{P}}U_{\Sigma }h\right) \left( \omega \right)  \nonumber \\
&=&\int_{\QTR{mathbb}{R}^{d}}\left( \int_{\QTR{mathbb}{R}}\left( h\left(
\omega +s\varepsilon _{x}\right) -h\left( \omega \right) \right) d\rho
\left( s\right) \right) \varphi \left( x\right) d\sigma \left( x\right) .
\label{2eq54}
\end{eqnarray}

\begin{example}
\begin{enumerate}
\item  Let $\rho =\varepsilon _{1}$, then $\mu _{\mathtt{CP}}=\pi _{\sigma }$
and, of course, (\ref{2eq54}) coincides with (\ref{2eq6}).

\item  Let $\rho =\frac{1}{2}\left( \varepsilon _{-1}+\varepsilon
_{1}\right) $ (for $d=1$ $\mu _{\mathtt{CP}}$ is generated by the so called
telegraph process) then the annihilation operator $\bigtriangledown
_{\varphi }^{\mathtt{CP}}$ has the form 
\[
\left( \bigtriangledown _{\varphi }^{\mathtt{CP}}h\right) \left( \omega
\right) =\int_{\QTR{mathbb}{R}}\left( \frac{1}{2}h\left( \omega +\varepsilon
_{x}\right) +\frac{1}{2}h\left( \omega -\varepsilon _{x}\right) -h\left(
\omega \right) \right) \varphi \left( x\right) d\sigma \left( x\right) . 
\]
\end{enumerate}
\end{example}

\begin{example}
Let $h\in \mathcal{F}C_{b}^{\infty }(\mathcal{D},\Omega )$ be given by 
\begin{eqnarray*}
h\left( \omega \right) =\exp \left( \left\langle \omega ,\log \left( 1+\eta
\right) \right\rangle -\left\langle \eta \right\rangle _{\sigma }\int_{%
\QTR{mathbb}{R}}sd\rho \left( s\right) \right) \\
=\exp \left( \left\langle \omega ,\log \left( 1+\eta \right) \right\rangle
-\left\langle \eta \right\rangle _{\sigma }m_{1}\left( \rho \right) \right)
\end{eqnarray*}
for $\mathcal{D}\ni \eta >-1.$ Then the annihilation operator $%
\bigtriangledown _{\varphi }^{\mathtt{CP}}$ applied to $h$ can be computed
to be 
\begin{eqnarray*}
\left( \bigtriangledown _{\varphi }^{\mathtt{CP}}h\right) \left( \omega
\right) &=&\left( \bigtriangledown _{\varphi }^{\mathtt{P}}\widehat{h}%
\right) \left( \widehat{\gamma }\right) \\
&=&\int_{\QTR{mathbb}{R}^{d+1}}\left( \widehat{h}\left( \widehat{\gamma }%
+\varepsilon _{\left( s,x\right) }\right) -\widehat{h}\left( \widehat{\gamma 
}\right) \right) \varphi \left( x\right) d\widehat{\sigma }\left( s,x\right)
\\
&=&\int_{\QTR{mathbb}{R}^{d+1}}\left( h\left( \omega +s\varepsilon
_{x}\right) -h\left( \omega \right) \right) \varphi \left( x\right) d%
\widehat{\sigma }\left( s,x\right) \\
&=&h\left( \omega \right) \int_{\QTR{mathbb}{R}^{d+1}}\left( \left( 1+\eta
\left( x\right) \right) ^{s}-1\right) \varphi \left( x\right) d\widehat{%
\sigma }\left( s,x\right) \\
&=&\left\langle \left( \left( 1+\eta \right) ^{\bullet }-1\right) ,\varphi
\right\rangle _{\widehat{\sigma }}h\left( \omega \right) .
\end{eqnarray*}
\end{example}

\begin{example}
Let $h\in \mathcal{F}C_{b}^{\infty }(\mathcal{D},\Omega )$ be given by 
\[
h(\omega )=\exp (-\left\langle \phi \right\rangle _{\rho }\left\langle
\varphi \right\rangle _{\widehat{\sigma }})\prod_{x\in \gamma }(1+\phi
(s_{x})\varphi (x)),\;\omega =\sum_{x\in \gamma }s_{x}\varepsilon _{x}\in
\Omega . 
\]
It is clear that $(U_{\Sigma }^{-1}e_{\pi _{\widehat{\sigma }}}^{\alpha }(%
\widehat{\varphi };\cdot ))(\omega )=h(\omega )$ $(\widehat{\varphi }:=\phi
\varphi )$ as it can be easily seen from the definitions of $U_{\Sigma
}^{-1} $ and $e_{\pi _{\widehat{\sigma }}}(\widehat{\varphi };\cdot )$ given
in (\ref{2eq71}) and (\ref{2eq2}), respectively. On the other hand from (\ref
{2eq44}) we have 
\[
(\nabla _{\widehat{\psi }}^{\mathtt{P}}e_{\pi _{\widehat{\sigma }}}^{\alpha
}(\widehat{\varphi };\cdot ))\left( \widehat{\gamma }\right) =(\widehat{%
\varphi },\widehat{\psi })_{L^{2}(\widehat{\sigma })}e_{\pi _{\widehat{%
\sigma }}}^{\alpha }(\widehat{\varphi };\widehat{\gamma }) 
\]
and therefore $(U_{\Sigma }^{-1}\bigtriangledown _{\widehat{\psi }}^{\mathtt{%
P}}U_{\Sigma }h)(\omega )=(\widehat{\varphi },\widehat{\psi })_{L^{2}(%
\widehat{\sigma })}h(\omega )$ which says that 
\[
(\bigtriangledown _{\psi }^{\mathtt{CP}}h)(\omega )=(\varphi ,\psi
)_{L^{2}(\sigma )}h(\omega ). 
\]
\end{example}

Now we proceed to compute an expression for the creation operator on
compound Poisson space.

\begin{proposition}
Let $g\in L^{2}(\Omega ,\mu _{\mathtt{CP}})$ be such that $U_{\Sigma }g\in 
\mathrm{Dom}(I_{\widehat{\sigma }}a^{+}(\widehat{\varphi })I_{\widehat{%
\sigma }}^{-1})$ and $\widehat{\varphi }\in \mathcal{D}(\QTR{mathbb}{R}%
^{d+1}).$ Then the operator $U_{\Sigma }^{-1}(\bigtriangledown _{\widehat{%
\varphi }}^{\mathtt{P}})^{*}U_{\Sigma }$ has the following representation 
\[
\left( U_{\Sigma }^{-1}(\bigtriangledown _{\widehat{\varphi }}^{\mathtt{P}%
})^{*}U_{\Sigma }g\right) \left( \omega \right) =\int_{\QTR{mathbb}{R}%
^{d+1}}g\left( \omega -s\varepsilon _{x}\right) \widehat{\varphi }\left(
s,x\right) d\widehat{\gamma }\left( s,x\right) -g\left( \omega \right)
\left\langle \widehat{\varphi }\right\rangle _{\widehat{\sigma }}. 
\]
\end{proposition}

\noindent \textbf{Proof.\ }We know from (\ref{2eq4}) that for any $\widehat{g%
}\in \mathrm{Dom}(I_{\widehat{\sigma }}a^{+}(\widehat{\varphi })I_{\widehat{%
\sigma }}^{-1})$ the creation operator $(\bigtriangledown _{\widehat{\varphi 
}}^{\mathtt{P}})^{*}$ on Poisson space $L^{2}(\widehat{\Gamma },\pi _{%
\widehat{\sigma }})$ has the form 
\[
\left( \left( \bigtriangledown _{\widehat{\varphi }}^{\mathtt{P}}\right) ^{*}%
\widehat{g}\right) \left( \widehat{\gamma }\right) =\int_{\QTR{mathbb}{R}%
^{d+1}}\widehat{g}\left( \widehat{\gamma }-\varepsilon _{\widehat{x}}\right) 
\widehat{\varphi }\left( \widehat{x}\right) d\widehat{\gamma }\left( 
\widehat{x}\right) -\widehat{g}\left( \widehat{\gamma }\right) \left\langle 
\widehat{\varphi }\right\rangle _{\widehat{\sigma }}. 
\]
On the other hand, 
\begin{eqnarray*}
\left( U_{\Sigma }^{-1}(\bigtriangledown _{\widehat{\varphi }}^{\mathtt{P}%
})^{*}U_{\Sigma }g\right) \left( \omega \right) &=&\left( (\bigtriangledown
_{\widehat{\varphi }}^{\mathtt{P}})^{*}\widehat{g}\right) \left( \widehat{%
\gamma }\right) \\
&=&\int_{\QTR{mathbb}{R}^{d+1}}\widehat{g}\left( \widehat{\gamma }%
-\varepsilon _{\widehat{x}}\right) \widehat{\varphi }\left( \widehat{x}%
\right) d\widehat{\gamma }\left( \widehat{x}\right) -\widehat{g}\left( 
\widehat{\gamma }\right) \left\langle \widehat{\varphi }\right\rangle _{%
\widehat{\sigma }} \\
&=&\int_{\QTR{mathbb}{R}^{d+1}}g\left( \omega -s\varepsilon _{x}\right) 
\widehat{\varphi }\left( s,x\right) d\widehat{\gamma }\left( s,x\right)
-g\left( \omega \right) \left\langle \widehat{\varphi }\right\rangle _{%
\widehat{\sigma }},
\end{eqnarray*}
which proves the result of the proposition.\hfill $\blacksquare $

As before if we choose $\widehat{\varphi }=1\varphi ,$ in the case when $%
\rho $ is finite, then we have the following form for the \textbf{creation
operator} $(\bigtriangledown _{\varphi }^{\mathtt{CP}})^{*}$ in compound
Poisson space $L^{2}(\Omega ,\mu _{\mathtt{CP}})$: 
\begin{eqnarray*}
\left( \left( \bigtriangledown _{\varphi }^{\mathtt{CP}}\right) ^{*}g\right)
\left( \omega \right) &\mbox{$:=$}&\left( U_{\Sigma }^{-1}\left(
\bigtriangledown _{\varphi }^{\mathtt{P}}\right) ^{*}U_{\Sigma }g\right)
\left( \omega \right) \\
&=&\int_{\QTR{mathbb}{R}^{d+1}}g\left( \omega -s\varepsilon _{x}\right)
\varphi \left( x\right) d\widehat{\gamma }\left( s,x\right) -g\left( \omega
\right) \rho (\QTR{mathbb}{R})\left\langle \varphi \right\rangle _{\sigma }.
\end{eqnarray*}

\begin{remark}
The generalized Charlier polynomials in $L^{2}(\pi _{\widehat{\sigma }}),$
according to (\ref{2eq56}), have the following representation 
\[
\left( \left( \nabla _{\widehat{\varphi }}^{\mathtt{P}}\right) ^{*n}1\right)
\left( \widehat{\gamma }\right) =\left\langle C_{n}^{\widehat{\sigma }%
}\left( \widehat{\gamma }\right) ,\widehat{\varphi }^{\otimes
n}\right\rangle . 
\]
Their images under $U_{\Sigma }^{-1}$ have the following form 
\begin{eqnarray*}
\left( U_{\Sigma }^{-1}\left\langle C_{n}^{\widehat{\sigma }}\left( \cdot
\right) ,\widehat{\varphi }^{\otimes n}\right\rangle \right) \left( \omega
\right) &=&\left\langle C_{n}^{\widehat{\sigma }}\left( \Sigma ^{-1}\omega
\right) ,\widehat{\varphi }^{\otimes n}\right\rangle \\
&=&\left( U_{\Sigma }^{-1}\left( \nabla _{\widehat{\varphi }}^{\mathtt{P}%
}\right) ^{*n}U_{\Sigma }1\right) \left( \omega \right) .
\end{eqnarray*}
In particular for finite measure $\rho $ and $\widehat{\varphi }=\varphi $
we obtain 
\[
\left( \left( \nabla _{\varphi }^{\mathtt{CP}}\right) ^{*n}1\right) \left(
\omega \right) =\left\langle C_{n}^{\widehat{\sigma }}\left( \Sigma
^{-1}\omega \right) ,\varphi ^{\otimes n}\right\rangle . 
\]
\end{remark}

\section{Gamma analysis\label{2eq66}}

\subsection{Definition and properties}

In this section we consider the classical (real) Schwartz triple 
\[
\mathcal{D}\left( \QTR{mathbb}{R}^{d}\right) =:\mathcal{D}\subset
L^{2}\left( \QTR{mathbb}{R}^{d}\right) \subset \mathcal{D}^{\prime }:=%
\mathcal{D}^{\prime }\left( \QTR{mathbb}{R}^{d}\right) .
\]

\begin{definition}
We call \textbf{\emph{Gamma noise} }the measure $\mu _{\mathtt{G}}^{\sigma
}\,$on the measure space $(\mathcal{D}^{\prime }{,}\mathcal{B}(\mathcal{D}%
^{\prime }))$ determined via its Laplace transform 
\begin{eqnarray*}
l_{\mu _{\mathtt{G}}^{\sigma }}\left( \varphi \right)  &=&\int_{\mathcal{D}%
^{\prime }}\exp \left( \left\langle \omega ,\varphi \right\rangle \right)
d\mu _{\mathtt{G}}^{\sigma }\left( \omega \right)  \\
&=&\exp \left( -\left\langle \log \left( 1-\varphi \right) \right\rangle
_{\sigma }\right) ,\;1>\varphi \in \mathcal{D}.
\end{eqnarray*}
\end{definition}

\begin{remark}
\label{2eq21}In order to apply the Minlos' theorem we note that $\mu _{%
\mathtt{G}}^{\sigma }\,$is a special case of $\mu _{\mathtt{CP}}$ for the
choice of $\rho $ as follows 
\begin{equation}
\rho \left( \Delta \right) =\int_{\Delta \cap \left] 0,\infty \right[ }\frac{%
e^{-s}}{s}ds,\;\Delta \in \mathcal{B}\left( \QTR{mathbb}{R}\right) .
\label{2eq67}
\end{equation}
Whence by Minlos' theorem $\mu _{\mathtt{G}}^{\sigma }$ is well-defined, of
course $l_{\mu _{\mathtt{G}}^{\sigma }}$ is an analytic function.
\end{remark}

\begin{remark}
\label{2eq68}Let us explain the term ``Gamma noise''. If $d=1$ and $\sigma =m$%
, then for any $t>0$ the value of the Laplace transform 
\[
l_{\mu _{\mathtt{G}}^{m}}\left( \lambda 1\!\!1_{\left[ 0,t\right] }\right)
=\exp \left( -t\log \left( 1-\lambda \right) \right) ,\;\lambda <1
\]
coincides with the Laplace transform $l_{\xi \left( t\right) }\left( \lambda
\right) $ of a random variable $\xi \left( t\right) $ having two-side Gamma
distribution, i.e., the density of the distribution function has the form 
\[
p_{t}\left( x\right) =\frac{1}{2}\frac{\left| x\right| ^{t-1}e^{-\left|
x\right| }}{\Gamma \left( t\right) },\;t>0.
\]
The process $\{\xi \left( t\right) ,t>0;\,\xi \left( 0\right) :=0\}$ is
known as Gamma process, see e.g. \cite[Section 19]{Ta67}. Thus the triple $(%
\mathcal{D}^{\prime },\mathcal{B}(\mathcal{D}^{\prime }),\mu _{\mathtt{G}%
}^{m})$ is a direct representation of the generalized stochastic process $\{%
\dot{\xi}\left( t\right) ,t\geq 0\}$ (detailed information on generalized
stochastic process can be found in \cite{GV68}) which is a distributional
derivative of the Gamma process $\{\xi \left( t\right) ,t\geq 0\}$. In other
words, the image of $\mu _{\mathtt{G}}^{m}$ under the transformation 
\[
\mathcal{D}^{\prime }\ni \omega \longmapsto G_{t}(\omega ):=\left\langle
\omega ,1\!\!1_{[0,t]}\right\rangle \in \QTR{mathbb}{R},\;t\in \QTR{mathbb}{R%
}_{+}
\]
coincides with the two-sided Gamma distribution, i.e., 
\[
(\mu _{\mathtt{G}}^{m}\circ G_{t}^{-1})(\Delta )=\int_{\Delta
}p_{t}(x)dm(x),\;\Delta \in \mathcal{B}(\QTR{mathbb}{R}).
\]
So the term ``Gamma noise'' is natural for $\mu _{\mathtt{G}}^{m}.$
\end{remark}

\begin{remark}
As $\mu _{\mathtt{G}}^{\sigma }$ is a special case of compound Poisson
measure one can obtain the representations of generalized Charlier
polynomials, annihilation and creation operators etc. along the lines of
Subsection ~\ref{2eq65}. It is worth noting that $\rho (\QTR{mathbb}{R}%
)=\infty $ nevertheless one can set in (\ref{2eq55}) $\widehat{\varphi }(%
\widehat{x})=1\varphi (x)$ and obtain the representation (\ref{2eq52}) 
\[
(\nabla _{\varphi }^{\mathtt{G}}h)(\omega )=\int_{\QTR{mathbb}{R}%
^{d}}\int_{0}^{\infty }(h(\omega +s\varepsilon _{x})-h(\omega ))\frac{e^{-s}%
}{s}ds\varphi (x)d\sigma (x) 
\]
for the annihilation operator in Gamma space. Indeed, by (\ref{2eq41}) 
\begin{eqnarray*}
h(\omega +s\varepsilon _{x})-h(\omega ) &=&H((\left\langle \omega ,\varphi
_{1}\right\rangle ,\ldots ,\left\langle \omega ,\varphi _{N}\right\rangle
)+s(\varphi _{1}(x),\ldots ,\varphi _{N}(x))) \\
&&-H(\left\langle \omega ,\varphi _{1}\right\rangle ,\ldots ,\left\langle
\omega ,\varphi _{N}\right\rangle ),\;H\in C_{b}^{\infty }(\QTR{mathbb}{R}%
^{N})
\end{eqnarray*}
whence by Lagrange theorem it follows that $|h(\omega +s\varepsilon
_{x})-h(\omega )|\leq Cs.$ Therefore the integral over $[0,\infty )$
converges and the r.h.s. of above equality is well-defined.
\end{remark}

\begin{remark}
Let us assume that $d=1$ and $\sigma =m.$ Then $\mu _{\mathtt{G}}^{m}$
corresponds to a distributional derivative of the Gamma process $\xi =\{\xi
(t),t\geq 0\}$ on a probability space $(\Omega ,\mathcal{F},P)$ (see Remark ~%
\ref{2eq68}). The Gamma process is L\'{e}vy one, such that $\QTR{mathbb}{E}%
[\xi (t)]=t$ and $\QTR{mathbb}{E}[(\xi (t)-t)^{2}]=t$, $t\geq 0.$ Thus the
centered Gamma process $\{\xi (t)-t,t\geq 0\}$ is a normal martingale and
one can define the $n$-multiple stochastic integrals $I_{n}(f^{(n)},\xi )$, $%
f^{(n)}\in \mathrm{Exp}_{n}L^{2}(m)$ with respect to $\xi $ and the space of
chaos decomposable random variables from $L^{2}(\Omega ,\mathcal{F},P)$: 
\[
\QTR{mathfrak}{C}(\xi )=\left\{ \sum_{n=0}^{\infty }I_{n}(f^{(n)},\xi
)\,|\,f^{(n)}\in \mathrm{Exp}_{n}L^{2}(m),\sum_{n=0}^{\infty }n!\left|
f^{(n)}\right| ^{2}<\infty \right\} 
\]
(for more details see \cite{M93}).

It follows from results of \cite{De90} that $\xi $ does not possess CRP,
i.e., $\QTR{mathfrak}{C}(\xi )$ is a proper subspace of $L^{2}(\Omega ,%
\mathcal{F}_{\xi },P)=:L^{2}(\xi )$, where $\mathcal{F}_{\xi }$ denotes a $%
\sigma $-algebra which is generated by the collection $\{\xi (t),t\geq 0\}$.
\end{remark}

\subsection{Chaos decomposition of Gamma space}

Let us now consider a function $\alpha :\mathcal{D}\longrightarrow \mathcal{D%
}$ defined by 
\[
\alpha \left( \varphi \right) \left( x\right) =\frac{\varphi \left( x\right) 
}{\varphi \left( x\right) -1},\;\varphi \in \mathcal{D},\,x\in \QTR{mathbb}{R%
}^{d}. 
\]
We stress that $\alpha $ is a holomorphic function on a neighborhood of zero 
$\mathcal{U}_{\alpha }\subset \mathcal{D},$ in other words $\alpha \in 
\mathrm{Hol}_{0}(\mathcal{D},\mathcal{D}).$

Because of the holomorphy of $l_{\mu _{\mathtt{G}}^{\sigma }}\,$and $l_{\mu
_{\mathtt{G}}^{\sigma }}\left( 0\right) =1,$ there exists a neighborhood of
zero $\mathcal{U}_{\alpha }^{\prime }\subset \mathcal{U}_{\alpha }$ such
that the normalized exponential $e_{\mu _{\mathtt{G}}^{\sigma }}^{\alpha
}(\varphi ;\omega )$ is holomorphic for any $\varphi \in \mathcal{U}_{\alpha
}^{\prime }$ and $\omega \in \mathcal{D}^{\prime }.$ Then 
\begin{eqnarray}
e_{\mu _{\mathtt{G}}^{\sigma }}^{\alpha }\left( \varphi ;\omega \right) &%
\mbox{$:=$}&\,\frac{\exp \left( \left\langle \omega ,\alpha \left( \varphi
\right) \right\rangle \right) }{l_{\mu _{\mathtt{G}}^{\sigma }}\left( \alpha
\left( \varphi \right) \right) }  \nonumber \\
&=&\exp \left( \left\langle \omega ,\frac{\varphi }{\varphi -1}\right\rangle
-\left\langle \log \left( 1-\varphi \right) \right\rangle _{\sigma }\right)
,\,\varphi \in \mathcal{U}_{\alpha }^{\prime }.  \label{2eq12}
\end{eqnarray}
We use the holomorphy of $\varphi \mapsto e_{\mu _{\mathtt{G}}^{\sigma
}}^{\alpha }(\varphi ;\omega )$ to expand it in a power series which, with
Cauchy's inequality, polarization identity and kernel theorem, give us 
\begin{equation}
e_{\mu _{\mathtt{G}}^{\sigma }}^{\alpha }\left( \varphi ;\omega \right)
=\sum_{n=0}^{\infty }\frac{1}{n!}\left\langle P_{n}^{\mu _{\mathtt{G}%
}^{\sigma },\alpha }\left( \omega \right) ,\varphi ^{\otimes n}\right\rangle
,\;\varphi \in \mathcal{U}_{\alpha }^{\prime },\,\omega \in \mathcal{D}%
^{\prime },  \label{2eq60}
\end{equation}
where $P_{n}^{\mu _{\mathtt{G}}^{\sigma },\alpha }:\mathcal{D}^{\prime
}\rightarrow \mathcal{D}^{\prime \widehat{\otimes }n}.$ $\{P_{n}^{\mu _{%
\mathtt{G}}^{\sigma },\alpha }\left( \cdot \right) =:L_{n}^{\sigma }\left(
\cdot \right) |n\in \QTR{mathbb}{N}_{0}\}$ is called the system of \textbf{%
generalized Laguerre kernels} on Gamma space $(\mathcal{D}^{\prime },%
\mathcal{B}(\mathcal{D}^{\prime }),\mu _{\mathtt{G}}^{\sigma }).$ From (\ref
{2eq60}) it follows immediately that for any $\varphi ^{(n)}\in \mathcal{D}^{%
\widehat{\otimes }n},$ $n\in \QTR{mathbb}{N}_{0}$ the function 
\[
\mathcal{D}^{\prime }\ni \omega \mapsto \left\langle L_{n}^{\sigma }\left(
\omega \right) ,\varphi ^{(n)}\right\rangle 
\]
is a polynomial of the order $n$ on $\mathcal{D}^{\prime }.$ The system of
functions 
\[
\left\{ L_{n}^{\sigma }\left( \varphi ^{(n)}\right) \left( \omega \right)
:=\left\langle L_{n}^{\sigma }\left( \omega \right) ,\varphi ^{\left(
n\right) }\right\rangle ,\,\forall \varphi ^{(n)}\in \mathcal{D}^{\widehat{%
\otimes }n},\,n\in \QTR{mathbb}{N}_{0}\right\} 
\]
is called the system of \textbf{generalized Laguerre polynomials} for the
Gamma measure $\mu _{\mathtt{G}}^{\sigma }$. In one-dimensional case this
system coincides with the system of even continuations of Laguerre
polynomials $\{L_{n}^{(t-1)},n\in \QTR{mathbb}{N}_{0}\}$ which are
orthogonal with respect to the density $p_{t}(x)$, see e.g. \cite{R71}.
Namely in the notation of Remark \ref{2eq68} the following equality holds 
\[
\left\langle L_{n}^{m}(\omega ),1\!\!1_{[0,t]}^{\otimes n}\right\rangle
=L_{n}^{(t-1)}(G_{t}(\omega )). 
\]

Now we proceed establishing the following result. Let $\varphi ,\psi \in 
\mathcal{U}_{\alpha }^{\prime }$ be given, then using (\ref{2eq12}) follows
that for $\lambda _{1},\lambda _{2}\in \QTR{mathbb}{R}$ 
\begin{eqnarray}
&&\int_{\mathcal{D}^{\prime }}e_{\mu _{\mathtt{G}}^{\sigma }}^{\alpha
}\left( \lambda _{1}\varphi ;\omega \right) e_{\mu _{\mathtt{G}}^{\sigma
}}^{\alpha }\left( \lambda _{2}\psi ;\omega \right) d\mu _{\mathtt{G}%
}^{\sigma }\left( \omega \right)  \nonumber \\
&=&\exp \left( \left\langle -\log \left( 1-\lambda _{1}\varphi \right) -\log
\left( 1-\lambda _{2}\psi \right) \right\rangle _{\sigma }\right)  \nonumber
\\
&&\cdot \int_{\mathcal{D}^{\prime }}\exp \left( \left\langle \omega ,\tfrac{%
\lambda _{1}\varphi }{\lambda _{1}\varphi -1}+\tfrac{\lambda _{2}\psi }{%
\lambda _{2}\psi -1}\right\rangle \right) d\mu _{\mathtt{G}}^{\sigma }\left(
\omega \right)  \nonumber \\
&=&\exp \left( \left\langle -\log \left( 1-\lambda _{1}\varphi \right) -\log
\left( 1-\lambda _{2}\psi \right) \right\rangle _{\sigma }\right)  \nonumber
\\
&&\cdot \exp \left( -\left\langle \log \left( 1-\tfrac{\lambda _{1}\varphi }{%
\lambda _{1}\varphi -1}-\tfrac{\lambda _{2}\psi }{\lambda _{2}\psi -1}%
\right) \right\rangle _{\sigma }\right)  \nonumber \\
&=&\exp \left( -\left\langle \log \left( 1-\lambda _{1}\varphi \lambda
_{2}\psi \right) \right\rangle _{\sigma }\right)  \nonumber \\
&=&l_{\mu _{\mathtt{G}}^{\sigma }}\left( \lambda _{1}\lambda _{2}\varphi
\psi \right) .  \label{2eq29}
\end{eqnarray}
Since $l_{\mu _{\mathtt{G}}^{\sigma }}\in \mathcal{M}_{a}(\mathcal{D}%
^{\prime }),$ then (\ref{2eq29}) turns out to be an analytic function on $%
\lambda _{1}$ and $\lambda _{2}$, hence 
\begin{equation}
l_{\mu _{\mathtt{G}}^{\sigma }}\left( \lambda _{1}\lambda _{2}\varphi \psi
\right) =\sum_{n=0}^{\infty }\frac{1}{n!}\left( \lambda _{1}\lambda
_{2}\right) ^{n}\left( \varphi ^{\otimes n},\psi ^{\otimes n}\right) _{%
\mathrm{Exp}_{n}^{\mathtt{G}}L^{2}\left( \sigma \right) },  \label{2eq13}
\end{equation}
where the coefficients $(\varphi ^{\otimes n},\psi ^{\otimes n})_{\mathrm{Exp%
}_{n}^{\mathtt{G}}L^{2}\left( \sigma \right) }$ are given by 
\[
\left( \varphi ^{\otimes n},\psi ^{\otimes n}\right) _{\mathrm{Exp}_{n}^{%
\mathtt{G}}L^{2}\left( \sigma \right) }=\left. \frac{d^{n}}{dt^{n}}\exp
\left( -\left\langle \log \left( 1-t\varphi \psi \right) \right\rangle
_{\sigma }\right) \right| _{t=0} 
\]
and $\mathrm{Exp}_{n}^{\mathtt{G}}L^{2}\left( \sigma \right) $ stands for a
quasi-$n$-particle subspace of $\mathrm{Exp}^{\mathtt{G}}L^{2}(\sigma )$
defined by (\ref{2eq58}) below.

By using the formula, see e.g. \cite{B58} and \cite{GR81}, 
\begin{eqnarray*}
&\,&\frac{d^{n}}{dt^{n}}e^{f\left( t\right) } \\
&=&\sum_{{\QATOP{i_{1}+2i_{2}+\cdots +ki_{k}=n}{i_{1},i_{2},\ldots ,i_{k}\in 
\QTR{mathbb}{N}_{0}}}}\frac{n!}{i_{1}!\ldots i_{k}!}\left( \frac{f^{\left(
1\right) }\left( t\right) }{1!}\right) ^{i_{1}}\left( \frac{f^{\left(
2\right) }\left( t\right) }{2!}\right) ^{i_{2}}\cdots \left( \frac{f^{\left(
k\right) }\left( t\right) }{k!}\right) ^{i_{k}}e^{f\left( t\right) }
\end{eqnarray*}

follows 
\begin{eqnarray}
&&\left( \varphi ^{\otimes n},\psi ^{\otimes n}\right) _{\mathrm{Exp}_{n}^{%
\mathtt{G}}L^{2}\left( \sigma \right) }  \nonumber \\
&=&\sum_{{\QATOP{i_{1}+2i_{2}+\cdots +ki_{k}=n}{i_{1},i_{2},\ldots ,i_{k}\in 
\QTR{mathbb}{N}_{0}}}}\frac{n!}{i_{1}!i_{2}!\ldots i_{k}!}\frac{1}{%
2^{i_{2}}\ldots k^{i_{k}}}  \nonumber \\
&&\cdot \left( \int_{\QTR{mathbb}{R}^{d}}\varphi \left( x\right) \psi \left(
x\right) d\sigma \left( x\right) \right) ^{i_{1}}\left( \int_{\QTR{mathbb}{R}%
^{d}}\varphi ^{2}\left( x\right) \psi ^{2}\left( x\right) d\sigma \left(
x\right) \right) ^{i_{2}}  \nonumber \\
&&\cdot \cdots \left( \int_{\QTR{mathbb}{R}^{d}}\varphi ^{k}\left( x\right)
\psi ^{k}\left( x\right) d\sigma \left( x\right) \right) ^{i_{k}}.
\label{2eq25}
\end{eqnarray}

On the other hand 
\begin{eqnarray}
&&\int_{\mathcal{D}^{\prime }}e_{\mu _{\mathtt{G}}^{\sigma }}^{\alpha
}\left( \lambda _{1}\varphi ;\omega \right) e_{\mu _{\mathtt{G}}^{\sigma
}}^{\alpha }\left( \lambda _{2}\psi ;\omega \right) d\mu _{\mathtt{G}%
}^{\sigma }\left( \omega \right)  \nonumber \\
&=&\sum_{n,m=0}^{\infty }\frac{\lambda _{1}^{n}\lambda _{2}^{m}}{n!m!}\int_{%
\mathcal{D}^{\prime }}\left\langle L_{n}^{\sigma }\left( \omega \right)
,\varphi ^{\otimes n}\right\rangle \left\langle L_{m}^{\sigma }\left( \omega
\right) ,\psi ^{\otimes m}\right\rangle d\mu _{\mathtt{G}}^{\sigma }\left(
\omega \right) .  \label{2eq14}
\end{eqnarray}
Then a comparison of coefficients between (\ref{2eq13}) and (\ref{2eq14})
gives us 
\[
\int_{\mathcal{D}^{\prime }}\left\langle L_{n}^{\sigma }\left( \omega
\right) ,\varphi ^{\otimes n}\right\rangle \left\langle L_{m}^{\sigma
}\left( \omega \right) ,\psi ^{\otimes m}\right\rangle d\mu _{\mathtt{G}%
}^{\sigma }\left( \omega \right) =\delta _{nm}n!\left( \varphi ^{\otimes
n},\psi ^{\otimes n}\right) _{\mathrm{Exp}_{n}^{\mathtt{G}}L^{2}\left(
\sigma \right) }, 
\]
which shows the orthogonality property of the system $\{L_{n}^{\sigma
}\left( \cdot \right) |n\in \QTR{mathbb}{N}_{0}\}.$

Since $(\cdot ,\cdot )_{\mathrm{Exp}_{n}^{\mathtt{G}}L^{2}\left( \sigma
\right) }$ is $n$-linear we can extend it by polarization, linearity and
continuity to general smooth kernels $\varphi ^{\left( n\right) },$ $\psi
^{\left( n\right) }\in \mathrm{Exp}_{n}^{\mathtt{G}}L^{2}\left( \sigma
\right) .$ To this end we proceed as follows.

First we consider a partition of the numbers $I_{n}:=\{1,2,\ldots ,n\}$ in 
\[
I_{n}=\bigcup_{\alpha }I_{\alpha }=:\mathcal{I}^{\left( n\right) }. 
\]
Then for each such partition $\mathcal{I}^{\left( n\right) },$ we define $%
i_{k}$ by 
\[
i_{k}:=\#\left\{ I_{\alpha }|\,\left| I_{\alpha }\right| =k\right\} ,\;1\leq
k\leq n. 
\]
Finally we define the contraction of the kernel $\varphi ^{\left( n\right) }$
w.r.t. $\mathcal{I}^{\left( n\right) }$ as 
\[
\varphi _{\mathcal{I}^{\left( n\right) }}^{\left( n\right) }\left(
x_{1},x_{2},\ldots ,x_{n}\right) :=\varphi ^{\left( n\right) }\left(
x_{i_{1}},x_{i_{2}},\ldots ,x_{i_{k}}\right) , 
\]
where $x_{i_{m}}=\left( x_{m},x_{m},\ldots ,x_{m}\right) $ $m$-times, $1\leq
m\leq k$.

Hence the inner product is given by 
\begin{eqnarray}
\hspace{-0.7cm} &&\left( \varphi ^{\left( n\right) },\psi ^{\left( n\right)
}\right) _{\mathrm{Exp}_{n}^{\mathtt{G}}L^{2}\left( \sigma \right) } 
\nonumber \\
\hspace{-0.7cm} &=&\sum_{\mathcal{I}^{\left( n\right) }}n!\left(
\prod_{k=1}^{n}\frac{1}{i_{k}!k^{i_{k}}}\right) \left( n!\prod_{k=1}^{n}%
\frac{1}{\left( k!\right) ^{i_{k}}i_{k}!}\right) ^{-1}  \nonumber \\
\hspace{-0.7cm} &&\cdot \int_{\QTR{mathbb}{R}^{dn}}\varphi _{\mathcal{I}%
^{\left( n\right) }}^{\left( n\right) }\left( x_{1},\ldots ,x_{n}\right)
\psi _{\mathcal{I}^{\left( n\right) }}^{\left( n\right) }\left( x_{1},\ldots
,x_{n}\right) d\sigma ^{\otimes n}\left( \vec{x}\right)  \nonumber \\
\hspace{-0.7cm} &=&\sum_{\mathcal{I}^{\left( n\right)
}}\prod_{k=1}^{n}((k-1)!)^{i_{k}}\int_{\QTR{mathbb}{R}^{dn}}\varphi _{%
\mathcal{I}^{\left( n\right) }}^{\left( n\right) }\left( x_{1},\ldots
,x_{n}\right) \psi _{\mathcal{I}^{\left( n\right) }}^{\left( n\right)
}\left( x_{1},\ldots ,x_{n}\right) d\sigma ^{\otimes n}\left( \vec{x}\right)
,  \label{2eq49}
\end{eqnarray}
where the sum extends over all possible partition $\mathcal{I}^{\left(
n\right) }$ of $I_{n}.$

Hence we have established the proposition.

\begin{proposition}
Let $\varphi ,\psi \in \mathcal{D}$ be given. Then the system of generalized
Laguerre polynomials verifies the following orthogonality property 
\[
\int_{\mathcal{D}^{\prime }}\left\langle L_{n}^{\sigma }\left( \omega
\right) ,\varphi ^{\left( n\right) }\right\rangle \left\langle L_{m}^{\sigma
}\left( \omega \right) ,\psi ^{\left( m\right) }\right\rangle d\mu _{\mathtt{%
G}}^{\sigma }\left( \omega \right) =\delta _{nm}n!\left( \varphi ^{\left(
n\right) },\psi ^{\left( n\right) }\right) _{\mathrm{Exp}_{n}^{\mathtt{G}%
}L^{2}\left( \sigma \right) }, 
\]
where $(\varphi ^{\left( n\right) },\psi ^{\left( n\right) })_{\mathrm{Exp}%
_{n}^{\mathtt{G}}L^{2}\left( \sigma \right) }$ is defined by (\ref{2eq49})
above.
\end{proposition}

As a consequence of the last proposition we have established the following
isomorphism 
\begin{equation}
I:L^{2}\left( \mu _{\mathtt{G}}^{\sigma }\right) \longrightarrow
\bigoplus_{n=0}^{\infty }\mathrm{Exp}_{n}^{\mathtt{G}}L^{2}\left( \sigma
\right) =:\mathrm{Exp}^{\mathtt{G}}L^{2}\left( \sigma \right) .  \label{2eq58}
\end{equation}
Therefore for any $F\in L^{2}(\mu _{\mathtt{G}}^{\sigma })$ there is a
sequence $(f^{\left( n\right) })_{n=0}^{\infty }\in \mathrm{Exp}^{\mathtt{G}%
}L^{2}(\sigma )$ such that 
\[
F\left( \omega \right) =\sum_{n=0}^{\infty }\left\langle L_{n}^{\sigma
}\left( \omega \right) ,f^{\left( n\right) }\right\rangle , 
\]
moreover 
\[
\left\| F\right\| _{L^{2}\left( \mu _{\mathtt{G}}^{\sigma }\right)
}=\sum_{n=0}^{\infty }n!\left| f^{\left( n\right) }\right| _{\mathrm{Exp}%
_{n}^{\mathtt{G}}L^{2}\left( \sigma \right) }^{2}. 
\]

\begin{remark}
Hence we see that the Gamma noise does not produce the standard Fock type
isomorphism since the inner product $\left( \cdot ,\cdot \right) _{\mathrm{%
Exp}_{n}^{\mathtt{G}}L^{2}\left( \sigma \right) }$ do not coincide with the
inner product in the $n$-particle subspace, $L^{2}(\sigma )^{\widehat{%
\otimes }n}$.
\end{remark}

\begin{remark}
The orthogonal polynomials of independent-increment processes (in
particular, Gamma-process) were constructed in \cite{KS76}. It is worth
noting that these polynomials of Gamma process differ from generalized
Laguerre polynomials.
\end{remark}

\noindent \textbf{Acknowledgments.} We are grateful to Dr.~E.~W.~Lytvynov
for helpful discussions. Financial support of the DFG through the project AL
214/9-2, DLR through WTZ-projects X271.7, JNICT - plurianual I\&D N.~219/94,
DAAD, and TMR Nr.~ERB4001GT957046 are gratefully acknowledged.

\end{document}